\documentclass[11pt,amsfonts]{article}
\usepackage{graphicx}
\usepackage{latexsym}
\usepackage{amssymb}
\usepackage{amsmath}
\usepackage{enumerate}
\usepackage{color}
\usepackage{layout}
\usepackage{eufrak}
\usepackage{float}
\usepackage{dsfont}

\newtheorem{prop}{Proposition}
\newtheorem{lemma}{Lemma}
\newtheorem{definition}{Definition}
\newtheorem{corollary}{Corollary}
\newtheorem{theorem}{Theorem}
\newtheorem{remark}{Remark}

\def\real{{\mathord{{\rm I\kern-2.8pt R}}}}        
\def\inte{{\mathord{{\rm I\kern-2.8pt N}}}}

\def\sZZ{{\rm Z\kern-2.8ptem{}Z}}

\def\z{{\mathchoice
  {\sZZ}
  {\sZZ}
  {\rm Z\kern-0.30em{}Z}
  {\rm Z\kern-0.25em{}Z} }}
\def\sQQ{{\kern 0.27em \vrule height1.45ex width0.03em depth0em
          \kern-0.30em \rm Q}}
\def\qu{{\mathchoice
    {\sQQ}
    {\sQQ}
  {\kern 0.225em \vrule height1.05ex width0.025em depth0em \kern-0.25em \rm Q}
  {\kern 0.180em \vrule height0.78ex width0.020em depth0em \kern-0.20em \rm Q}
        }}
\def\sCC{{\kern 0.27em \vrule height1.45ex width0.03em depth0em
          \kern-0.30em \rm C}}
\def\complex{{\mathchoice
    {\sCC}
    {\sCC}
  {\kern 0.225em \vrule height1.05ex width0.025em depth0em \kern-0.25em \rm C}
  {\kern 0.180em \vrule height0.78ex width0.020em depth0em \kern-0.20em \rm C}
        }}

\newcommand{\ba}{\begin{array}}
\newcommand{\ea}{\end{array}}
\newcommand{\be}{\begin{equation}}
\newcommand{\ee}{\end{equation}}
\newcommand{\bea}{\begin{eqnarray}}
\newcommand{\eea}{\end{eqnarray}}
\newcommand{\beaa}{\begin{eqnarray*}}
\newcommand{\eeaa}{\end{eqnarray*}}

\def\b{\beta}

\def\z{\zeta}

\font\tenmath=msbm10 \font\sevenmath=msbm7 \font\fivemath=msbm5
\newfam\mathfam \textfont\mathfam=\tenmath
\scriptfont\mathfam=\sevenmath \scriptscriptfont\mathfam=\fivemath

\def \b{\noindent}

\def \={{\buildrel {\rm (law)} \over =}}

\def\qed{ \hfill \vrule width.25cm height.25cm depth0cm\smallskip}

\newcommand{\basa}{\begin{assumption}}
\newcommand{\easa}{\end{assumption}}

\newcommand{\bas}{\begin{assum}}
\newcommand{\eas}{\end{assum}}

\def\limsup{\mathop{\overline{\rm lim}}}

\newcommand{\ignore}[1]{}
\begin{document}

\renewcommand{\thefootnote}{\fnsymbol{footnote}}

\renewcommand{\thefootnote}{\fnsymbol{footnote}}

\title{Generalized $k$-variations and Hurst parameter estimation for the fractional wave equation via Malliavin calculus }
\author{Radomyra Shevchenko $^{\dagger}$, Meryem Slaoui $^{*}$ and C. A. Tudor  $^{*}$\vspace*{0.2in} \\
 $^{\dagger}$ Fakult\"at f\"ur Mathematik,
LSIV, TU Dortmund \\
 Vogelpothsweg 87, 44227 Dortmund, Germany. \\
\quad radomyra.shevchenko@tu-dortmund.de \vspace*{0.1in}\\
$^{*}$ Laboratoire Paul Painlev\'e, Universit\'e de Lille \\
 F-59655 Villeneuve d'Ascq, France.\\
\quad meryem.slaoui@univ-lille.fr\\
 \quad ciprian.tudor@math.univ-lille.fr\\
\vspace*{0.1in} }

\maketitle

\begin{abstract} 
We analyze {\it  the generalized $k$-variations} for the solution to the wave equation driven by an additive Gaussian noise which behaves as a fractional Brownian with Hurst parameter $H>\frac{1}{2}$ in time and which  is white in space. The $k$-variations are defined along {\it filters} of any order $p\geq 1$ and of any length. We show that the sequence of generalized $k$-variation satisfies a Central Limit Theorem when $p> H+\frac{1}{4}$ and we estimate the rate of convergence for it via the Stein-Malliavin calculus. The results are applied to the estimation of the Hurst index. We construct several  consistent estimators for $H$ and these estimators are analyzed theoretically and numerically.
\end{abstract}

\vskip0.3cm

{\bf 2010 AMS Classification Numbers:}  60G15, 60H05, 60G18.

\vskip0.3cm

{\bf Key Words and Phrases}: Hurst parameter estimation; Fractional Brownian motion; Stochastic wave equation; Generalized
variation;  Stein-Malliavin calculus; Central
limit theorem.

\section{Introduction}
Since several decades, the statistical inference in stochastic (partial) differential equations (S(P)DE in the sequel) constitutes an intensive research direction in probability theory and mathematical statistics. Nowadays, a particular case  of wide interest is represented by the S(P)DE driven by fractional Brownian motion (fBm) and related processes, due to the vast area of application of such stochastic models. Many recent  works concern the estimation of the drift parameter for stochastic equation driven by fractional Brownian motion (we refer, among many others, to \cite{AV}, \cite{HN}, \cite{KB}, \cite{TV}), while fewer works deal with the estimation of the Hurst parameter in such stochastic equations.

In this paper, we consider the stochastic wave equation driven by  an additive Gaussian noise which behaves as a fractional Brownian motion in time and as a Wiener process  in space (we call it {\it fractional-white noise}). Our purpose is to construct and analyze an estimator for the Hurst parameter of the solution to this SPDE based on the observation of the solution at a fixed time and at a discrete number of points in space. The wave equation with fractional noise in time and/or in space has been studied in several works, such as \cite{BCA}, \cite{DT}, \cite{GLT}, \cite{KT}, \cite{QT} etc. We will use a standard method to construct estimators for the Hurst parameter, which is based on the $k$-variations of the observed process. The method has been recently employed in \cite{KT} for the case of quadratic variations, i.e. $k=2$. As for the fBm,  it was shown that the standard quadratic variation estimator is not asymptotically normal when the Hurst index becomes bigger than $\frac{3}{4}$ and this is inconvenient for statistical applications. In order to avoid this restriction and to get an estimator which is asymptotically Gaussian for every $H$, we will use the generalized $k$-variations, which basically means that the usual increment of the process is replaced by a higher order increment. The idea comes from the reference \cite{istas} and since it has been used by many authors (see e.g. \cite{coeur} or \cite{CTV}). More precisely, if $(u(t,x), t\geq 0, x\in \mathbb{R})$ denotes the solution to the wave equation with fractional-white noise, we define the (centered)  generalized $k$-variation statistics ($k\geq 1$ integer),
\begin{equation}\label{vn-intro}
V_{N} (k, \alpha) = \frac{1}{N-l} \sum_{i=l} ^{N} \left[ \frac{ \left|  U ^{\alpha} \left(\frac{i}{N} \right) \right| ^{k} }{\mathbf{E}\left|  U ^{\alpha} \left(\frac{i}{N} \right) \right| ^{k}}-1\right],
\end{equation}
where $U ^{\alpha} \left(\frac{i}{N}\right)$ represents the spatial increment of the solution $u$ at $\frac{i}{N}$ along a {\it filter} $\alpha$ of power (order) $p\geq 1$ and length $ l+1\geq 1$ (see the next section for the precise definition). 

By using chaos expansion and the recent Stein-Malliavin calculus we show that the sequence (\ref{vn-intro}) satisfies a Central Limit Theorem (CLT) as $N\to \infty$ (in the spirit of \cite{Breuer}) whenever $p> H+\frac{1}{4}$ and in this way the restriction $H<\frac{3}{4}
$  can be avoided by choosing a filter of order $p\geq 2$,  i. e.  by replacing, for example, the usual increment by a higher order increment. We will obtain the rate of convergence under the Wasserstein distance for this convergence in law and we also prove a multidimensional CLT. So we generalize the findings in \cite{KT} to filters of any power $p\geq 1$ and to $k$-variations of any order $k\geq 1$ and in addition we show that in the special case $p=1 $ and $  H>\frac{3}{4}$ a non-Gaussian limit theorem occurs with limit distribution related to the Rosenblatt distribution (but more complex than it). 

These theoretical results are then applied to the estimation of Hurst index of the solution the the fractional-white wave equation. Based on the behavior of the sequence  (\ref{vn-intro}), we prove that the associated $k$-variation estimator for $H$ is consistent and asymptotically normal. Moreover, we provide a numerical analysis of the estimator when $k=2$ by analysing its performance on various filters and for several values of the Hurst parameter and confirming via simulation the theoretical results.

We organized the paper as follows. Section 2 contains some preliminaries. We present in this part the basic facts concerning the solution to the fractional-white wave equation, we introduce the filters and the increment of the solution along filters. In Section 3, we prove a CLT for the sequence (\ref{vn-intro}) for any integer $k\geq 1$ and we obtain the rate of convergence when $k$ is even via the Stein-Malliavin theory. In Section 4, we show a non-central limit limit in the case $k=2, H>\frac{3}{4}$ and for filters of order $p=1$.  Section 5 concerns the estimation of the Hurst parameter of the solution to the fractional-white wave equation (\ref{systeme wave}). We included here theoretical results related to the behavior of the $k$-variations estimators for the Hurst index as well as simulations and numerical analysis for the performance of the estimators. Section 6 (the Appendix) contains the basic tools from Malliavin calculus needed in the paper.

\section{Preliminaries}

We introduce here the fractional-white heat equation and its solution and we present the basic definitions and the notation concerning the filters used in our work.

\subsection{The solution to the wave equation with fractional-colored noise}

The object of our study will be the solution to the following stochastic wave equation
\begin{equation}
\left\{
\begin{array}{rcl}\label{systeme wave}
\frac{\partial^2 u}{\partial t^2}(t,x)&=&\Delta
u(t,x)+\dot W^H(t,x),\;t\geq 0,\;x \in \mathbb{R}^{d},\,d\geq 1,\nonumber\\
\noalign{\vskip 2mm}
u(0, x)&=& 0, \quad x \in \mathbb{R}^{d},\nonumber\\
\noalign{\vskip 2mm} \frac{\partial u}{\partial t}(0,x) &=& 0,\quad
x \in \mathbb{R} ^{d},
\end{array} \right.
\end{equation}
where  $\Delta$ is the Laplacian on $\mathbb{R}^{d}$, $d \geqslant 1$, and $W^{H}$ is a fractional-white Gaussian noise which is defined as a real valued centered Gaussian field $W_{H}=\{ {W_{t}^{H} (A); t \in  [0, T ], A \in B_{b} (\mathbb{R}^{d} )}\}$, over a given complete filtered probability space
$(\Omega, \mathcal{F}, ( \mathcal{F}_{t} )_{t\geqslant0} , \mathbb{P})$, with covariance function given by 
\begin{equation}\label{covW}
\mathbf{E}\left(W_{t}^{H} (A)W_{s}^ {H} (B)\right) = R_{H} (t, s)\lambda(A \cap B), \forall A, B \in \mathcal{B}_{b} (\mathbb{R}^{d} ),
\end{equation}
where $R_{H}$ is the covariance of the fractional brownian motion
 \begin{equation*}
R_{H}(t,s)=\frac{1}{2}\left(t^{2H}+ s^{2H}- \vert t-s \vert^{2H} \right), \hskip0.3cm s,t\geq 0.
\end{equation*}
We denoted by $B_{b} (\mathbb{R}^{d}) $ the class of bounded Borel subsets of $\mathbb{R}^{d}$ and we will assume throughout this work $H\in \left(\frac{1}{2}, 1\right).$

The solution of the equation (\ref{systeme wave})  is understood in the mild sense, that is, it is defined as   a square-integrable
centered field $u=\left( u(t,x);\;t\in[0,T], x \in\mathbb{R}^d \right)$
 defined by
\begin{equation}\label{sol-wave-WFrac}
u(t,x)=\int_{0}^{t}\int_{\mathbb{R}^d}G_1(t-s,x-y)W^H(\mathrm{d}s,\mathrm{d}y),\quad t\geq 0, x\in \mathbb{R} ^ {d},
\end{equation}
where $G_{1}$ is the fundamental solution to the wave equation and the integral in (\ref{sol-wave-WFrac}) is a Wiener integral with respect to the Gaussian process $W^{H}$. Recall that for $d=1$ (we will later restrict to this situation in our work) we have, for $t\geq 0$ and $x\in \mathbb{R} $,
\begin{equation}
\label{g1}
G_1(t,x)=\frac{1}{2}\mathds{1}_{\{|x|<t\}}.
\end{equation}

 \b We know (see e.g. \cite{BCA}) that the solution (\ref{sol-wave-WFrac}) is well-defined if and only if
$$d< 2H+1$$
and it is self-similar in time and stationary in space. Other properties of the solution can be found in \cite{BCA}, \cite{DT} or \cite{T}. In particular, the spatial covariance of the solution can be expressed as follows
\begin{eqnarray}\label{cov-mild}
\mathbf{E}\left( u(t, x)u(t, y)\right)
&=&\frac{1}{2}\left(c_{H} \vert y -x\vert^{ 2H+1} - t\frac{\vert y - x\vert^{ 2H}}{2}+ \frac{t^{2H+1}}{2H+1}\right)1 _{\{\vert y -x\vert<t\}}\nonumber\\
&&+
\frac{(2t - \vert y -x\vert)^{2H+1}} {8(2H+1)}1 _{\{t\leq \vert y -x\vert<2t\}}
\end{eqnarray}
with $c_{H}=\frac{4H-1}{4(2H+1)}$. When $t>1$ and $x,y \in [0,1]$, this expression reduces to 
\begin{eqnarray}\label{cov-mild-simp}
\mathbf{E}\left( u(t, x)u(t, y)\right)
&=&\frac{1}{2}\left(c_{H} \vert y -x\vert^{ 2H+1} - t\frac{\vert y - x\vert^{ 2H}}{2}+ \frac{t^{2H+1}}{2H+1}\right) \text{ for }x,\, y\in [0,\,1].
\end{eqnarray}

We will fix for the rest of the work $t >1$ and we will associate to the process $\left( u(t,x), x\in [0,1]\right) $ its canonical Hilbert space $\mathcal{H}$ which is defined as the closure of the  linear space generated by the indicator functions $\{1_{[0, x]}, x\in [0,1]\}$ with respect to the inner product 
\begin{eqnarray*}
\langle1_{[0,x]},1_{[0,y]}\rangle_{\mathcal{H}}&=& \mathbf{E}\left(u(t,x)u(t,y)\right).
\end{eqnarray*}
We will denote by $I_{q}$ the multiple stochastic integral of order $q\geq 1$ with respect to the Gaussian process $\left( u(t,x), x\in [0,1]\right) $ and by $D$ the Malliavin derivative with respect to this process. We refer to the Appendix for the basic elements of the Malliavin calculus.

We will also use in Section \ref{nclt} multiple stochastic integrals with respect to the fractional-white noise $W^{H}$ with covariance (\ref{covW}). We  use the notation $ I ^{W}_{q}$ to indicate the multiple integral of order $q\geq 1$ with respect to $ W^{H}$. 

\subsection{Filters}
In this paragraph we will define the filters and the increments of the solution to (\ref{systeme wave}) along filters. We start with several definitions and notations needed along this paper.

\begin{definition} \label{def1}
Given $ l \in \mathbb{N}$ and $p \in \mathbb{N}^{*}$,  a vector $\alpha= (\alpha_{0}, . . . , \alpha_{l})$ is called a filter of length $l + 1\geq 1$ and order  (or power) $p\geq 1$ such that

\begin{equation*}
\left\{
  \begin{array}{rcr}
   \sum_{q=0}^{l}\alpha_{q} q^{r} & =  0,& 0\leq r \leq p-1, \\
   \sum_{q=0}^{l}\alpha_{q} q^{p} & \neq  0 \\
  \end{array}
\right.
\end{equation*}

with the convention $0^{0}=1$.
\end{definition}

For a filter $ \alpha =(a_{0}, a_{1},.., a_{l}) $ of length $l+1\geq 1$ and of order $p\geq 1$ we define the space-filtered process (or the spatial increment of the process $u$ along the filter $\alpha$)
\begin{equation}\label{uin}
U ^ {\alpha} \left( \frac{i}{N} \right) = \sum_{r=0} ^ {l} a_{r} u \left( t, \frac{i-r}{N} \right) \mbox{ for } i=l,.., N.
\end{equation}

\b We denote for $j\geq 1$
\begin{equation*}
\pi ^ {\alpha, N }_{H}(j):= \mathbf{E} \left[ U ^ {\alpha} \left( \frac{i}{N} \right)U ^ {\alpha} \left( \frac{i+j}{N} \right)\right].
\end{equation*}
From the covariance formula (\ref{cov-mild})  we can write
\begin{eqnarray}
\pi ^ {\alpha, N }_{H}(j)&=& \sum_{r_{1}, r_{2}=0} ^ {l} a_{r_{1}} a_{r_{2}} \mathbf{E} \left[ u\left(t, \frac{ i-r_{1}}{N}\right) u\left(t, \frac{ i+j-r_{2}}{N}\right) \right]\nonumber\\
&=&k_{1} \frac{ 1}{ N ^ {2H}} \Phi_{H, \alpha } (j)+ k_{2} \frac{1}{N ^ {2H+1}} \Phi_{H+\frac{1}{2}, \alpha } (j)\label{14d-1}
\end{eqnarray}
with
\begin{equation*}
\Phi_{H, \alpha} (j) =  \sum_{r_{1}, r_{2}=0} ^ {l} a_{r_{1}} a_{r_{2}}\vert j+r_{1}-r_{2} \vert ^ {2H}, \hskip0.5cm j\geq 0
\end{equation*}
and $k_{1}=-\frac{t}{4} $ and $k_{2}=\frac{c_{H}}{2}=\frac{4H-1}{8(2H+1)}$. We write for further use

\begin{eqnarray}\label{c1}
c_{1}(H)=\frac{-t}{4}\sum_{q,r=0}^{l}\alpha_{q}\alpha_{r} \vert q-r \vert^{2H},
\mbox{ and }
c_{2}(H)=\frac{c_{H}}{2}\sum_{q,r=0}^{l}\alpha_{q}\alpha_{r} \vert q-r \vert^{2H+1}.
\end{eqnarray}
In particular, from (\ref{14d-1})
\begin{eqnarray*}
\pi ^ {\alpha, N }_{H}(0)&=& \mathbf{E}\left[  U ^ {\alpha} \left( \frac{i}{N} \right)\right] ^ {2} = k_{1} \frac{ 1}{ N ^ {2H}} \Phi_{H, \alpha } (0)+ k_{2} \frac{1}{N ^ {2H+1}} \Phi_{H+\frac{1}{2}, \alpha } (0)\\
&=& c_{1}(H) \frac{ 1}{ N ^ {2H}} + c_{2}(H) \frac{ 1}{ N ^ {2H+1}}
\end{eqnarray*}

We will need the below  technical lemma to etablish the asymptotical equivalent of $\Phi_{H,\,\alpha}$ and similar expressions. The  proof of the lemma is based on a Taylor expansion, see \cite{coeur} or \cite{istas}.
\begin{lemma}\label{lemma}
Let $H \in \mathbb R^+\backslash \mathbb N$  and $\alpha^{(1)}$, $\alpha^{(2)}$ be filters of lengths $l_1+1$, $l_2+2$ and of orders $p_1,\,p_2 \geq 1 $ respectively. Then 

\begin{eqnarray*}
\sum_{q=0}^{l_1}\sum_{r=0}^{l_2}\alpha^{(1)}_{q}\alpha^{(2)}_{r}\vert q-r +k\vert ^{2H} \underset{k \to \infty}\sim & \kappa_{H} k^{2H-2p}
\end{eqnarray*}
with $\kappa_{H}=\sum_{q=0}^{l_1}\sum_{r=0}^{l_2}\alpha^{(1)}_{q}\alpha^{(2)}_{r}\frac{2H(2H-1) \ldots (2H-2p +1)}{2p !} (q-r)^{2p}$, where $p=\min (p_1,\,p_2)$.
\end{lemma}
In the sequel, we write  $a_{k}\sim _{k\to \infty} b_{k}$ to indicate  that the sequences $a_{k}, b_{k}$ have the same behavior as $k\to \infty$. 

\section{Central limit theorem for the spatial  $k$-variations}\label{clt}

In this section we focus on the asymptotic behavior in distribution of the $k$-variation in space of the solution  to the fractional-white wave equation, defined via a filter of power $p\geq 1$.  In the first step we show the $k$-variation satisfies a CLT when $p>H+\frac{1}{4}$. Next, by taking $k$ to be an even integer, we derive a Berry-Ess\'een type bound for this convergence in distribution via the Stein-Malliavin  calculus. Restricting ourselves in addition to $k=2$, we prove a multidimensional CLT, which is needed for the estimation of the Hurst parameter.  

\subsection{Central Limit Theorem}

Fix $t>1$ and let $\alpha$ be a filter of length $l+1\geq 1$ and of power $p\geq 1$ as in Definition \ref{def1}. Let $u$ be given by (\ref{sol-wave-WFrac}). For any integer $k\geq 1$ we define the centered spatial $k$-variations of the process $\left( u(t, x), x\in \mathbb{R} \right)$ by

\begin{equation}\label{vna}
V_{N} (k, \alpha) = \frac{1}{N-l} \sum_{i=l} ^{N} \left[ \frac{ \left|  U ^{\alpha} \left(\frac{i}{N} \right) \right| ^{k} }{\mathbf{E}\left|  U ^{\alpha} \left(\frac{i}{N} \right) \right| ^{k}}-1\right]
\end{equation}
with $U ^{\alpha} \left(\frac{i}{N} \right) $ given by (\ref{uin}). We will show that the sequence (\ref{vna}) satisfies a CLT.  In order to do this we will use a criterion based on Malliavin calculus. The first step is to expand in chaos the $k$-variation sequence $V_{N} (k, \alpha)$.  Noticing  that the filtered process $U^{\alpha}$ as a linear combination of centered Gaussian random variables is a centered Gaussian process, we get 
\begin{equation}\label{p1-1}
\mathbf{E}\left(U^{\alpha}\left( \frac{i}{N}\right)^{k}\right)=E_{k}\mathbf{E}\left(U^{\alpha}\left( \frac{i}{N}\right)^{2}\right)^{\frac{k}{2}},
\end{equation}
where $E_{k}$ denotes  the $k$-th absolute moment of a standard Gaussian variable given by $E_{k}=\frac{2^\frac{k}{2}\Gamma(\frac{k+1}{2})}{\Gamma(\frac{1}{2})}$.
We introduce the variable 
\begin{equation}\label{p1-2}
Z^{\alpha}\left( \frac{i}{N}\right)=\frac{U^{\alpha}\left( \frac{i}{N}\right)}{(\pi ^ {\alpha, N }_{H}(0))^{1\slash 2}}.
\end{equation}
It is clear that $Z^{\alpha}\left( \frac{i}{N}\right)$ is a standard Gaussian variable and $\operatorname{Corr}\left(Z^{\alpha}\left( \frac{i}{N}\right), Z^{\alpha}\left( \frac{j}{N}\right)\right)=\operatorname{Corr}\left(U^{\alpha}\left( \frac{i}{N}\right), U^{\alpha}\left( \frac{j}{N}\right)\right)$, where $\operatorname{Corr}$ denotes the correlation coefficient. 
Using (\ref{p1-1}) and (\ref{p1-2}) we can write $V_{N}$ as follows:

\begin{eqnarray*}
V_{N}(k,\alpha)&=&\frac{1}{N-l}\sum_{i=l}^{N}\left[\frac{\vert U^{\alpha}(\frac{i}{N})\vert^{k}}{\mathbf{E}\vert U^{\alpha}(\frac{i}{N})\vert^{k}}-1\right]=\frac{1}{N-l}\sum_{i=l}^{N}\left[\frac{\vert Z^{\alpha}(\frac{i}{N})\vert^{k}}{E_{k}}-1\right].\\
\end{eqnarray*} 

\b Recall the expansion of the development in Hermite polynomials of  the function $H^{k}(t)=\frac{\vert t \vert ^{k}}{E_{k}}-1$   given in Lemma 2 of \cite{coeur},
\begin{equation*}
H^{k}(t)= \sum_{j=1}^{\infty}c_{j}^kH_{j}(t),
\end{equation*}
where  $c_{2j+1}^k=0$ for $j \geqslant0$, $c^k_{2j}=\frac{1}{(2j)!}\prod_{i=0}^{j-1}(k-2i)$ for $j \geqslant 1$ and $H_{j}(t)$ denotes the j-th Hermite polynomial  defined by 
\begin{equation*}
H_{j}(t)=\sum_{a=0}^{[\frac{j}{2}]} (-1)^{a} \frac{a!}{(j-2a)!a!}2^{-a}t^{j-2a}.
\end{equation*}

\b Observing that for $$C_{i,\,\alpha}:=\sum_{q=0}^{l}\alpha_{q}\mathbf{1}_{\left[0,\frac{i-q}{N}\right]}$$ we have from (\ref{14d-1}) that $\left\|\frac{C_{i,\,\alpha}}{(\pi ^ {\alpha, N }_{H}(0))^{1\slash 2}} \right\|_{\mathcal H} = 1$ we can express $Z^{\alpha}\left(\frac{i}{N} \right)$  as an integral with respect to the process $(u(t,x), x \in \mathbb{R})$ since 
the increment $u(t,y)-u(t,x)$ can be
expressed as ${I_{1}}(\mathds{1}_{[x,y]})$ (recall that $I_{1}$  represents the multiple integral of order 1 with respect to the Gaussian process $\left(u(t,x), x\in [0,1]\right)$) for every $x<y$:\

\begin{eqnarray*}
Z^{\alpha}\left(\frac{i}{N} \right)&=&I_{1}\left(\frac{C_{i,\,\alpha}}{(\pi ^ {\alpha, N }_{H}(0))^{1\slash 2}}\right).\\
\end{eqnarray*}
Since we have $H_{q}(I_{1}(h))=\frac{1}{q!}I_{q}(h^{\otimes q}) $ for $\Vert h \Vert _{\mathcal{H}}=1$ we get

\begin{eqnarray*}
V_{N}(k,\alpha)&=&\frac{1}{N-l}\sum_{i=l}^{N}H^{k}\left(Z^{\alpha}\left(\frac{i}{N}\right)\right)
=\frac{1}{N-l}\sum_{q\geqslant 1 }c^k_{2q}\sum_{i=l}^{N}H_{2q}\left(Z^{\alpha}\left(\frac{i}{N}\right)\right)\\
&=&\frac{1}{N-l}\sum_{q\geqslant 1 }c^k_{2q}\sum_{i=l}^{N} H_{2q}\left( I_{1}\left(\frac{C_{i,\,\alpha}}{(\pi ^ {\alpha, N }_{H}(0))^{1\slash 2}}\right)\right)\\
&=&\frac{1}{N-l}\sum_{q\geqslant 1 }\frac{c^k_{2q}}{(2q)!}\sum_{i=l}^{N}I_{2q}\left( \left(\frac{C_{i,\,\alpha}}{(\pi ^ {\alpha, N }_{H}(0))^{1\slash 2}}\right)^{\otimes 2 q}\right).
\end{eqnarray*}
Hence, we obtain the following chaotic expansion of the $k$-variation statistics

\begin{equation}\label{vnka}
V_{N} (k, \alpha) = \frac{1}{N-l} \sum_{i=l} ^{N} \sum _{q=1} ^{\infty} \frac{ c_{2q} ^{k}}{(2q)! } I_{2q} \left(    \frac{ C_{i, \alpha} ^{\otimes 2q}} {( \pi ^{\alpha, N}_{H} (0) ) ^{q}}  \right)= \sum_{q\geq 1} I_{2q} (f_{N, 2q}) 
\end{equation}
with
\begin{equation}
\label{fnq}
f_{N, 2q} =\frac{ c_{2q} ^{k}}{(2q)! }  \frac{1}{N-l} \sum_{i=l} ^{N}    \frac{ C_{i, \alpha} ^{\otimes 2q}} {( \pi ^{\alpha, N}_{H} (0) ) ^{q}}.
\end{equation}

Let us start by analyzing the asymptotic behavior of the mean square of each kernel $f_{N, 2q}$.

\begin{lemma}\label{5d-4}
For $N, q\geq 1$, let $f_{N, 2q}$ be given by (\ref{fnq}). Then
\begin{equation*}
(N-l) (2q)! \Vert f_{N, 2q}\Vert ^ {2}_{ \mathcal{H} ^ {\otimes 2q } } \to _{N \to \infty} \frac{ (c_{2q} ^{k}) ^ {2} }{(2q)! }  \sum_{v \in \mathbb{Z}}( \varphi_{H, \alpha} (v) ) ^ {2q}:= \sigma _{2q} ^ {2} 
\end{equation*}
for $H<p-\frac{1}{4q}$ (i.e. $H < 1- \frac{1}{4q}$ for $p = 1$ and $H \in (\frac{1}{2},1)$ for $p\geq2$ ), where we use the notation
\begin{equation}\label{fi}
 \varphi _{H, \alpha } (v)=\frac{ \Phi_{H, \alpha }(v)}{\Phi_{H, \alpha} (0) }. 
\end{equation}
Moreover, $\sigma ^ {2}:=\sum_{q\geq 1} \sigma _{2q} ^ {2} <\infty. $ For $p=q=1$, $H=3\slash 4$
\begin{equation}\label{14d-3}
\frac{N-l}{\log (N-l)}2! \Vert f_{N, 2}\Vert ^ {2}_{ \mathcal{H} ^ {\otimes 2 } } \to _{N \to \infty} c^2:= \frac{(c_{2}^{k}) ^{2}} {2}\lim_{N} \log N \sum_{\vert v\vert \leq N} (\rho_{H, \alpha} (v) )^{2} <\infty .
\end{equation}
\end{lemma}
{\bf Proof: } From (\ref{fnq}), we get
\begin{eqnarray*}
 (2q)! \Vert f_{N, 2q}\Vert ^ {2}_{ \mathcal{H} ^ {\otimes 2q } } &=& \frac{ (c_{2q} ^{k}) ^ {2} }{(2q)! }\frac{1}{ (N-l) ^ {2}}\sum_{i, j=l} ^ {N}  \frac{ \langle C_{i, \alpha}, C_{j, \alpha} \rangle _{\mathcal{H}} ^{2q}}{ ( \pi ^{\alpha, N}_{H} (0) ) ^{2q}} \\
&=& \frac{1}{(N-l) ^{2} }\frac{ (c_{2q} ^{k}) ^{2}}{(2q)! } \sum_{i,j=l} ^{N} \left( \rho^{\alpha, N}_{H} (j-i) \right) ^{2q},
\end{eqnarray*}
where we used the notation 
\begin{equation}
\label{ro}
\rho^{\alpha, N}_{H}(v)= \frac{ \pi ^{\alpha, N}_{H}(v) } { \pi ^{\alpha, N}_{H} (0)  } \mbox{ for } v\in \mathbb{Z}. 
\end{equation}
Next, we write

\begin{equation*}
\frac{1}{N-l} \sum_{i,j=l} ^{N} \left( \rho^{\alpha, N}_{H} (j-i) \right) ^{2q}=\sum _{v\in \mathbb{Z} } \left( \rho ^{\alpha, N}_{H} (v)\right) ^{2q} 1_{\{ \vert v\vert \leq N-l \}} \frac{N-\vert v\vert -l}{N-l}
\end{equation*}
and thus
\begin{eqnarray}\label{5d-1}
 (N-l)(2q)! \Vert f_{N, 2q}\Vert ^ {2}_{ \mathcal{H} ^ {\otimes 2q } } &=& \frac{ (c_{2q} ^{k}) ^{2}}{(2q)! } \sum_{v \in \mathbb{Z}}  \left( \rho^{\alpha, N}_{H} (v) \right) ^{2q}1_{\{ \vert v\vert \leq N-l \}} \frac{ N-\vert v\vert -l}{N-l}.
\end{eqnarray}
Using the expression
\begin{eqnarray*}
 \rho ^{\alpha, N}_{H} (v)= \frac{ k_{1} \Phi _{H, \alpha} (v)N ^{-2H}+ k_{2} \Phi _{H+\frac{1}{2}, \alpha} (v) N^{-2H-1}}{ k_{1} \Phi _{H, \alpha} (0)N ^{-2H}+ k_{2} \Phi _{H+\frac{1}{2}, \alpha} (0) N^{-2H-1}}
= \frac{ \Phi_{H, \alpha} (v)+ a_{N}(v)}{\Phi_{H, \alpha} (0)+ a_{N}(0)}
\end{eqnarray*}
with
\begin{equation}\label{24n-2}
a_{N}(v) =\frac{k_{2}}{k_{1}N}  \Phi_{H+\frac{1}{2}, \alpha} (v) 
\end{equation}
we can write, with $\varphi_{H, \alpha}$ and $\rho^{\alpha, N}_{H}$ given by (\ref{fi}) and (\ref{ro}) respectively,
\begin{equation}\label{6d-1}
 b_{N, H} (v):= \rho ^{\alpha, N}_{H} (v) - \varphi_{H,\,\alpha} (v)
\end{equation}
and remark that due to Lemma \ref{lemma} for $v$ large enough
\begin{equation}
\label{bn}
\vert b_{N, H}(v) \vert \underset{\vert v\vert  \to \infty} \sim \left| a_{N}(v) \frac{1}{\Phi_{H,\,\alpha} (0)+a_N(0)}\right| \leq C \frac{1}{N} v ^ { 2H+1-2p},
\end{equation}
where $C>0$ does not depend on $N, v$.With this notation we can write
\begin{eqnarray*}
 (N-l)(2q)! \Vert f_{N, 2q}\Vert ^ {2}_{ \mathcal{H} ^ {\otimes 2q } } &=& \frac{ (c_{2q} ^{k}) ^{2}}{(2q)! } \sum_{v \in \mathbb{Z}}  \left( \varphi_{H, \alpha} (v)+ b_{N, H} (v)\right) ^{2q}1_{ \{\vert v\vert \leq N-l \}}  \frac{ N-\vert v\vert -l}{N-l}\\
&=& \frac{ (c_{2q} ^{k}) ^{2}}{(2q)! } \sum_{v \in \mathbb{Z}} \sum_{m=0} ^ {2q} \binom{2q}{m} \varphi_{H, \alpha }(v) ^ {m} (b_{N, H}(v)) ^ {2q-m}1_{ \{\vert v\vert \leq N-l \}}  \frac{ N-\vert v\vert -l}{N-l}\\
&=& \frac{ (c_{2q} ^{k}) ^{2}}{(2q)! } \sum_{v \in \mathbb{Z}}  \varphi_{H, \alpha}(v) ^ {2q} 1_{ \{\vert v\vert \leq N-l \}}  \frac{ N-\vert v\vert -l}{N-l}+ r_{N, q, 1} 
\end{eqnarray*}
with
\begin{equation}\label{rn1}
r_{N, q, 1}=  \frac{ (c_{2q} ^{k}) ^{2}}{(2q)! } \sum_{v \in \mathbb{Z}} \sum_{m=0} ^ {2q-1} \binom{2q}{m} \varphi_{H, \alpha }(v) ^ {m} (b_{N, H}(v)) ^ {2q-m}1_{ \{\vert v\vert \leq N-l \}}  \frac{ N-\vert v\vert -l}{N-l}.
\end{equation}
Clearly, by the dominated convergence theorem
$$\frac{ (c_{2q} ^{k}) ^{2}}{(2q)! } \sum_{v \in \mathbb{Z}}  \varphi_{H, \alpha}(v) ^ {2q} 1_{ \{\vert v\vert \leq N-l \}}  \frac{ N-\vert v\vert -l}{N-l}\to _{N\to \infty}\sigma _{2q} ^ {2},$$
which by Lemma \ref{lemma} is finite if $p=1$, $H<1-\frac{1}{4q}$, and for all $H\in (1\slash 2,\,1)$ in the other cases.
For $q=p=1$, $H=3\slash 4$,
$$ \frac{1}{\log (N-l)}\sum_{v \in \mathbb{Z}}  \varphi_{H, \alpha}(v) ^ {2} 1_{ \{\vert v\vert \leq N-l \}}  \frac{ N-\vert v\vert -l}{N-l}$$
converges to a positive constant and thus (\ref{14d-3}) is obtained.

In order to conclude it remains to show that the rest term $r_{N, q, 1}$ (\ref{rn1}) converges to $0$ as $N\to \infty$, for every $q\geq 1$. From (\ref{rn1}), using the bound (\ref{bn}) and Lemma \ref{lemma}, we have the estimate
\begin{eqnarray*}
\vert r_{N, q, 1}\vert  &\leq & C \sum_{m=0} ^ {2q-1} \binom{2q}{m}\frac{1}{N^ {2q-m}} \sum_{1\leq v\leq N-l} \vert v\vert  ^ {(2H-2) m } \vert v \vert ^ {(2H+1-2p) (2q-m)}:= \sum_{m=0} ^ {2q-1} r_{N,q, 1, m}
\end{eqnarray*}
and for each $m=0,.., 2q-1$, 
\begin{equation*}
r_{N, q, 1,m} \leq \frac{C}{N^{2q-m}}  \sum_{1\leq v\leq N-l}\vert v \vert ^ { (2H-2p) 2q +2q -m}. 
\end{equation*}
When the series $\sum_{v \in \mathbb{Z}} \vert v \vert ^ { (2H-2p) 2q +2q -m}$ converges we get 
$$r_{N, q, 1, m} \leq C \frac{1}{N ^ {2q-m}} \leq C\frac{1}{N} \to _{N\to \infty}0$$
and when the series diverges, 
$$r_{N,q,1, m} \leq C \frac{1}{N^ {2q-m}}  N ^ {(2H-2p)2q+2q-m+1} \leq C N ^ {(2H-2p) 2q+1} \to _{N\to \infty}0$$
if  $p=1, H < 1-\frac{1}{4q}$, or
$ p \geq 2$ and $H \in (\frac{1}{2},1)$. If $p=q=1$, $H=3\slash 4$, the quantity
$$\frac{1}{\log (N-l)} r_{N, q, 1, m}$$
will also converge to zero for $m=0,1$ using again (\ref{bn}) and Lemma \ref{lemma}.

The fact that the series $\sigma ^ {2}=\sum_{q\geq 1} \sigma_{2q} ^ {2} <\infty$ for $H<p-\frac{1}{4q}$ follows from the study of the $k$-variations of the fractional Brownian motion, see \cite{coeur} or \cite{N-P}.
\qed

\vskip0.2cm

We will consider the renormalized $k$-variation sequence
\begin{equation}
\label{gn}
G_{N} (k, \alpha)= \sqrt{N-l} V_{N} (k, \alpha).
\end{equation}
From the above Lemma \ref{5d-4}, it follows that 
$$\mathbf{E}\left[ G_{N}(k, \alpha)\right]^2 \to _{N \to \infty} \sigma ^ {2}$$ with $\sigma ^ {2} $ given in the statement of Lemma \ref{5d-4}. We will show that the sequence (\ref{gn}) satisfies a central limit theorem.

\begin{theorem}For a filter $\alpha$ of order $p\geq 1$ and of length $l+1\geq 1$, with $p>H+\frac{1}{4}$, let $G_{N}(k, \alpha) $ be given by (\ref{gn}). Then the sequence $ \left( G_{N}(k, \alpha) \right) _{N \geq 1}$ converges in distribution, as $N\to \infty$, to the Gaussian law $N(0, \sigma ^ {2})$. Moreover, for $p=1$, $H=3\slash 4$, the sequence  $ \left( \frac{1}{\sqrt{\log (N-l)}}G_{N}(k, \alpha) \right) _{N \geq 1}$ converges in distribution to $N(0,\,c^2)$. The constants $\sigma ^ {2}, c^ {2}$ are those appearing in Lemma \ref{5d-4}.
\end{theorem}
{\bf Proof: } Notice that from (\ref{vnka}), we can write 
\begin{equation}
\label{5d-3}
G_{N} (k, \alpha)=\sum_{q\geq1} I_{2q} (g_{N, 2q} ) \mbox{ with } g_{N, 2q} =\sqrt{N-l} f_{N, 2q}
\end{equation} 
with $f_{N, 2q}$ given by (\ref{fnq}). Our main tool to prove the asymptotic normality of (\ref{5d-3}) is Theorem 6.3.1 from \cite{N-P}. According to it, for $p>H+1\slash 4$ it suffices to show that
\begin{enumerate}
\item $(2q)! \Vert g_{N, 2q} \Vert ^ {2}_{\mathcal{H} ^ {\otimes 2q}}\to _{N\to \infty} \sigma _{2q}^ {2} $ and $\sigma ^ {2}:= \sum _{q\geq 1} \sigma _{2q}^ {2}<\infty , $

\item for every $q\geq 1$ and $r=1,.., 2q-1$, $\Vert g_{N, 2q}\otimes _{r} g_{N, 2q} \Vert _{ \mathcal{H} ^ {\otimes 4q-2r}} \to _{N \to \infty}0,$ 

\item $\lim _{M \to \infty} \sup_{N\geq 1} \sum_{q\geq M+1} (2q)! \Vert g_{N, 2q}\Vert ^ {2} _{\mathcal{H} ^ {\otimes 2q}}=0$
\end{enumerate}
and for $p=1$, $H=3\slash 4$,
\begin{enumerate}
\item $\frac{1}{\log (N-l)}(2q)! \Vert g_{N, 2q} \Vert ^ {2}_{\mathcal{H} ^ {\otimes 2q}}\to _{N\to \infty} 1_{\{q=1\}}c^2$,

\item for every $q\geq 1$ and $r=1,.., 2q-1$, $\frac{1}{\log (N-l)} \Vert g_{N, 2q}\otimes _{r} g_{N, 2q} \Vert _{ \mathcal{H} ^ {\otimes 4q-2r}} \to _{N \to \infty}0,$ 

\item $\lim _{M \to \infty} \sup_{N\geq 1} \sum_{q\geq M+1}\frac{1}{\log (N-l)}  (2q)! \Vert g_{N, 2q}\Vert ^ {2} _{\mathcal{H} ^ {\otimes 2q}}=0$.
\end{enumerate}
Point 1. in both cases follows from Lemma \ref{5d-4}.  Let us check point 2. By the definition of the contraction (\ref{contra}), we have for $q\geq 1$ and $r=1,.., 2q-1$

\begin{eqnarray*}
g_{N, 2q} \otimes _{r} g_{N, 2q}& =& \frac{1}{N-l}   \frac{ (c_{2q} ^{k}) ^{2}}{(2q)! }\sum_{i,j=l} ^ {N} \frac{ \langle C_{i, \alpha} , C_{j, \alpha}\rangle _{\mathcal{H}} ^ {r}}{ \pi ^ {\alpha, N}_{H}(0) ^ {2q}}C_{i, \alpha} ^ {\otimes 2q-r}\otimes C_{j, \alpha} ^ {\otimes 2q-r} 
\end{eqnarray*}
and
\begin{eqnarray*}
&&\Vert g_{N, 2q} \otimes _{r} g_{N, 2q}\Vert ^ {2}_{\mathcal{H}^ {\otimes 4q-2r}}\\
&=&\left(  \frac{ (c_{2q} ^{k}) ^{2}}{(2q)! }\right) ^ {2} \frac{1}{(N-l) ^ {2}}\sum_{i_{1}, i_{2}, i_{3}, i_{4}=l}^ {N}\frac{  \langle C_{i_{1}, \alpha}, C_{i_{2}, \alpha}\rangle^ {2q-r}  _{\mathcal{H}}  \langle C_{i_{2}, \alpha}, C_{i_{3}, \alpha}\rangle^ {r}  _{\mathcal{H}}  \langle C_{i_{3}, \alpha}, C_{i_{4}, \alpha}\rangle^ {2q-r}  _{\mathcal{H}}  \langle C_{i_{4}, \alpha}, C_{i_{1}, \alpha}\rangle^ {r}  _{\mathcal{H}} }{ \pi ^ {\alpha, N}_{H}(0) ^ {2q}}\\
&=&\left(  \frac{ (c_{2q} ^{k}) ^{2}}{(2q)! }\right) ^ {2} \frac{1}{(N-l) ^ {2}}\sum_{i_{1}, i_{2}, i_{3}, i_{4}=l}^ {N}\rho_{H} ^ {\alpha, N} (i_{1}-i_{2})^ {2q-r}  \rho_{H} ^ {\alpha, N} (i_{2}-i_{3})^ {r} \rho_{H} ^ {\alpha, N} (i_{3}-i_{4})^ {2q-r} \rho_{H} ^ {\alpha, N} (i_{4}-i_{1})^ {r} 
\end{eqnarray*}
with $\rho_{H} ^ {\alpha, N}$ given by (\ref{ro}). We use the fact that
\begin{eqnarray*}
\sum_{i_{1}, i_{2}, i_{3}, i_{4}=l}^ {N}&& \rho_{H} ^ {\alpha, N} (i_{1}-i_{2})^ {2q-r}  \rho_{H} ^ {\alpha, N} (i_{2}-i_{3})^ {r} \rho_{H} ^ {\alpha, N} (i_{3}-i_{4})^ {2q-r} \rho_{H} ^ {\alpha, N} (i_{4}-i_{1})^ {r} \\
&\leq & \sum_{n,\,m = l}^N \left(\left( \rho_{H} ^ {\alpha, N}1_{\{|\cdotp|\leq N-l\}}\right)^{2q-r}\ast  \left( \rho_{H} ^ {\alpha, N}1_{\{|\cdotp|\leq N-l\}}\right)^{r}\right)^2 (n-m),
\end{eqnarray*}
(where $\ast$ denotes convolution on $\mathbb Z$) and we obtain
\begin{eqnarray*}
&&\Vert g_{N, 2q} \otimes _{r} g_{N, 2q}\Vert ^ {2}_{\mathcal{H}^ {\otimes 4q-2r}}\\
&\leq & C \frac{1}{N-l} \sum_{v = l}^N \left(\left( \rho_{H} ^ {\alpha, N}1_{\{|\cdotp|\leq N-l\}}\right)^{2q-r}\ast  \left( \rho_{H} ^ {\alpha, N}1_{\{|\cdotp|\leq N-l\}}\right)^{r}\right)^2 (v)\\
&\leq & C\frac{1}{N-l} \left\|\left( \rho_{H} ^ {\alpha, N}1_{\{|\cdotp|\leq N-l\}}\right)^{2q-r} \right\|^2_{ l ^{4\slash 3}(\mathbb Z)} \left\|\left( \rho_{H} ^ {\alpha, N}1_{\{|\cdotp|\leq N-l\}}\right)^{r} \right\|^2_{ l ^{4\slash 3}(\mathbb Z)}\\
&=& C\frac{1}{N-l} \left(\sum_{|v|\leq N-l}\left( \rho_{H} ^ {\alpha, N}(v)\right) ^{(2q-r)\frac{4}{3}}\right)^{3\slash 2} \left(\sum_{|v|\leq N-l} \left(\rho_{H} ^ {\alpha, N}(v)\right) ^{r\frac{4}{3}}\right)^{3\slash 2}
\end{eqnarray*}
by virtue of the Young's inequality as in \cite{KT}. Note that for $v$ large enough we have by virtue of (\ref{bn})
\[b_{N,\,H}(v)1_{\{|v|\leq N-l\}}\leq C \frac{1}{N} v^{2H+1-2p}1_{\{|v|\leq N-l\}}\leq C v^{2H-2p}\leq C\varphi _H(v)1_{\{|v|\leq N-l\}},\]
and since all the powers invloved above are positive, this allows us to replace $\rho_{H} ^ {\alpha, N}$ with $\varphi _H$. Thus, for large $N$ the  norm $\Vert g_{N, 2q} \otimes _{r} g_{N, 2q}\Vert ^ {2}_{\mathcal{H}^ {\otimes 4q-2r}}$ is bounded by
$$C\frac{1}{N-l} \left(\sum_{|v|\leq N-l}|v|^{(2H-2p)(2q-r)\frac{4}{3}}\right)^{3\slash 2} \left(\sum_{|v|\leq N-l} |v|^{(2H-2p)r\frac{4}{3}}\right)^{3\slash 2}.$$
For $p\geq 2$ all these series converge. For $p=1$ and $H\leq \frac{3}{4}$ the only cases in which some of the series do not converege are $r=2q-1$ and $r=1$. However, the observation
\[\frac{1}{N-l}\sum_{{|v|\leq N-l}}|v|^{(2H-2)\frac{4}{3}}\sum_{{|v|\leq N-l}}|v|^{(2H-2)\frac{4}{3}}\leq C N^{-1}N^{\frac{8}{3}H-\frac{5}{3}}N^{\frac{8}{3}H-\frac{5}{3}}\to 0\]
ensures that even in those cases the term  $\Vert g_{N, 2q} \otimes _{r} g_{N, 2q}\Vert ^ {2}_{\mathcal{H}^ {\otimes 4q-2r}}$ converges to zero.

Concerning 3., fix $M\geq 1$ and recall that from (\ref{5d-1})
$$ (2q)! \Vert g_{N, 2q}\Vert ^ {2}_{ \mathcal{H} ^ {\otimes 2q } } = \frac{ (c_{2q} ^{k}) ^{2}}{(2q)! } \sum_{v \in \mathbb{Z}}  \left( \rho^{\alpha, N}_{H} (v) \right) ^{2q}1_{\{\vert v\vert \leq N-l\}} \frac{ N-\vert v\vert -l}{N-l}$$
and therefore, since $\vert  \rho^{\alpha, N}_{H} (v) \vert \leq 1$ for $|v|$ large enough, 
\begin{eqnarray*}
\sup_{N\geq 1} \sum_{q\geq M+1} (2q)! \Vert g_{N, 2q}\Vert ^ {2}_{ \mathcal{H} ^ {\otimes 2q } }&\leq & \sup_{N\geq 1} \sum_{q\geq M+1} \frac{ (c_{2q} ^{k}) ^{2}}{(2q)! } \sum_{v \in \mathbb{Z}}  \left( \rho^{\alpha, N}_{H} (v) \right) ^{2}1_{\{\vert v\vert \leq N-l\}} \frac{ N-\vert v\vert -l}{N-l}\\
&\leq &C \sum_{q\geq M+1} \frac{ (c_{2q} ^{k}) ^{2}}{(2q)! } \sum_{v \in \mathbb{Z}}  \varphi_{H, \alpha}(v) ^ {2} \\
&&+ C \sup_{N\geq 1} \sum_{q\geq M+1} \frac{ (c_{2q} ^{k}) ^{2}}{(2q)! } \sum_{v \in \mathbb{Z}} b_{N, H} (v) ^ {2} 1_{\{\vert v\vert \leq N-l\}} \frac{ N-\vert v\vert -l}{N-l}.
\end{eqnarray*}
From (\ref{bn}) 
$$ b_{N, H} (v) ^ {2}\leq C \frac{1}{N^ {2}} \mbox{ if } p\geq 2 $$
and
$$ \sum_{|v|\leq N-l }b_{N, H} (v) ^ {2} \frac{ N-\vert v\vert -l}{N-l} \leq C\frac{1}{ N ^ {2}} \sum _{|v|\leq N-l} v ^ {(2H-1)2}\leq C N ^ {4H-3} \mbox{ if } p=1, H<\frac{3}{4}.$$
Consequently, 
$$\sup_{N\geq 1} \sum_{q\geq M+1} (2q)! \Vert g_{N, 2q}\Vert ^ {2}_{ \mathcal{H} ^ {\otimes 2q } } \leq C   \sum_{q\geq M+1}  \frac{ (c_{2q} ^{k}) ^{2}}{(2q)! } \sum_{v \in \mathbb{Z}}  \varphi_{H, \alpha}(v) ^ {2} $$ and this tends to zero as $M \to \infty$, due to  the  convergence of the series $\sum_{q\geq 1}  \frac{ (c_{2q} ^{k}) ^{2}}{(2q)! }$.

For $\frac{1}{\sqrt{\log (N-l)}}G_N(k,\,\alpha)$ there is nothing to show since the case $q=1$ does not contribute to the limit.
 \qed

\subsection{Rate of convergence for even power variations}
Let $k\geq 2$ be an even integer.
Consider the sequence $G_{N}(k,\,\alpha)$ defined by (\ref{gn}). From (\ref{vnka}), since the coefficients $c_{2j}^ {k}$ vanish if $2j>k$, we get 
\begin{equation}\label{fnka2}
G_{N} (k, \alpha) = \frac{1}{\sqrt{N-l}} \sum_{i=l} ^{N} \sum _{q=1} ^{\frac{k}{2}} \frac{ c_{2q} ^{k}}{(2q)! } I_{2q} \left(    \frac{ C_{i, \alpha} ^{\otimes 2q}} {( \pi ^{\alpha, N}_{H} (0) ) ^{q}}  \right).
\end{equation}
Denote for every $q=1,2,.., \frac{k}{2}$ the $2q$-th chaos component of $G_{N} (k, \alpha)$ by
\begin{equation}\label{gnq}
G_{N} ^ {(2q)}(k, \alpha)= \frac{1}{\sqrt{N-l}} \sum_{i=l} ^{N} \frac{ c_{2q} ^{k}}{(2q)! } I_{2q} \left(    \frac{ C_{i, \alpha} ^{\otimes 2q}} {( \pi ^{\alpha, N}_{H} (0) ) ^{q}}  \right)= I_{2q}(g_{N, 2q})
\end{equation}
with $g_{N, 2q}$ from (\ref{5d-3}). Let us  consider the $\frac{k}{2}$-dimensional  random vector
\begin{equation*}
\mathbf{G} _{N}(k,\,\alpha) := \left( G_{N} ^ {(2)}(k, \alpha), G_{N} ^ {(4)}(k, \alpha),...., G_{N} ^ {(k)}(k, \alpha)  \right).
\end{equation*}
Notice that for every $q_{1}, q_{2} =1,.., \frac{k}{2}$  with $q_{1}\not= q_{2}$
$$\mathbf{E} G_{N} ^ {(2q_{1})}(k, \alpha)G_{N} ^ {(2q_{2})}(k, \alpha)=0$$
while for $ q_{1}=q_{2}=q$

\begin{equation*}
\mathbf{E} \left[ G_{N} ^ {(2q)}(k, \alpha)  \right]^ {2} = \frac{ (c_{2q} ^{k}) ^{2}}{(2q)! }\sum_{v\in \mathbb{Z}}\rho^{\alpha ,\,N}_H (v) ^{2q} 1_{\{\vert v\vert \leq N-l\}}\left( 1-\frac{\vert v\vert}{N-l} \right).
\end{equation*}
Let us introduce the matrix $C=(C_{q_{1}, q_{2}})_{q_{1}, q_{2}= 1,.., \frac{k}{2}}$ with components $C_{q_{1}, q_{2}}=0$ if $q_{1}\not=q_{2} $ and 
\begin{equation}
\label{cqq}
 C_{q, q}= \frac{ (c_{2q} ^{k}) ^{2}}{(2q)! }\sum_{v\in \mathbb{Z}}\varphi_{H, \alpha} (v) ^{2q} .
\end{equation}

In order to obtain convergence rates in the above CLT in terms of the Wasserstein distance we will use Corollary 3.6 from \cite{N-P-R} to prove a CLT as well as a rate of convergence for the vector $\mathbf{G}_{N} (k,\,\alpha)$, which will provide corresponding results for the $k$-variation statistics $V_{N}(k,  \alpha) $. For the sake of completeness we cite this corollary here.

\begin{corollary}
Fix $d\geq 2$ and $1\leq q_1\leq\dots \leq q_d$. Consider a vector $F:=(F_1,\dots F_d)=(I_{q_1}(f_1)\dots I_{q_d}(f_d))$ with $f_i\in\mathcal H ^{\odot q_i}$ for any $i=1,\dots , d$. Let $Z\sim N_d (0,\,C)$ with $C$ positive definite. Then
\[d_W(F,\,Z)\leq c \sqrt{\sum_{1\leq i,\,j\leq d}\mathbb \mathbf{E} \left[\left(C_{ij}-\frac{1}{q_j}\langle DF_i,\, DF_j \rangle_{\mathcal H}\right)^2\right]}\]
for some constant strictly positive $c$.
\end{corollary}
(In the one-dimensional case for a standard normal $Z$ this result is also true and can be found in \cite{N-P}. For $k=2$ the norming condition is satisfied, and the corollary is applicable.)

In order to apply the corollary for $F_{i}= G_{N} ^{ (2i)}, i=1,.., \frac{k}{2}$,  we will write each summand as
\begin{small}
\begin{eqnarray}
&&\mathbf{E} \left[\left(C_{ij}-\frac{1}{q_j}\langle DF_i,\, DF_j \rangle_{\mathcal H}\right)^2\right]\nonumber \\
&\leq & 2 \left(C_{ij}-\frac{1}{q_j}\mathbf{E} [\langle DF_i,\, DF_j \rangle_{\mathcal H}]\right)^2 + 2 \mathbf{E}  \left[\left(\frac{1}{q_j} \mathbf{E} [\langle DF_i,\, DF_j \rangle_{\mathcal H}]-\frac{1}{q_j}\langle DF_i,\, DF_j \rangle_{\mathcal H}\right)^2\right]\label{15d-1}
\end{eqnarray}
\end{small}
and conduct separate calculations for both parts. We start with a lemma for the deterministic part.

\begin{lemma}\label{det}
Let  $ G_{N} ^{(2q)}, C_{q, q} $ be given by (\ref{gnq}), (\ref{cqq}) respectively and assume $p\geq 2$. For $N$ large enough and for every $q=1,.., \frac{k}{2}$ we have for every $H\in \left( \frac{1}{2}, 1\right)$, 
$$\left| \mathbf{E}\left[ G_{N} ^ {(2q)}(k, \alpha) ^ {2} \right] - C_{q,q} \right| \leq C\frac{1}{N}. $$
For $p=1$ we have for $H\in \left(\frac{1}{2}, \frac{3}{4}\right)$
$$\left| \mathbf{E} \left[ G_{N} ^ {(2q)}(k, \alpha) ^ {2} \right] - C_{q,q} \right| \leq CN^{4H-3}. $$
\end{lemma}
{\bf Proof: } As in the proof of Lemma \ref{5d-4}, we have
%
\begin{eqnarray}
\mathbf{E} \left[ G_{N} ^ {(2q)}(k, \alpha) ^ {2} \right]&=& \frac{ (c_{2q} ^{k}) ^{2}}{(2q)! } \sum _{v\in \mathbb{Z} } \left( \rho ^{\alpha, N}_{H} (v)\right) ^{2q} 1_{ \{ \vert v\vert \leq N-l\}} \frac{N-\vert v\vert  -l }{N-l}.\label{24n-1}
\end{eqnarray}
Recall the representation $\rho^{\alpha ,\, N}_H(v)=\varphi_H(v)+b_{N,\,H}(v)$ introduced in Lemma \ref{5d-4}. We obtain by the Newton's formula
\begin{eqnarray*}
\mathbf{E} \left[ G_{N} ^ {(2q)}(k, \alpha) ^ {2} \right]&=& \frac{ (c_{2q} ^{k}) ^{2}}{(2q)! }\sum_{m=0} ^{2q} \binom{2q}{m} \sum_{v\in \mathbb{Z}}\varphi_{H, \alpha }(v) ^{m} b_{N,\,H}(v) ^{2q-m} 1_{ \{\vert v\vert \leq N-l \}} \frac{N-\vert v\vert -l}{N-l}\\
&=& \frac{ (c_{2q} ^{k}) ^{2}}{(2q)! }\sum_{v  \in \mathbb{Z}} \varphi_{H, \alpha} (v) ^{2q} 1_{\{ \vert v\vert \leq N-l\}} \frac{N-\vert v\vert-l }{N-l}\\
&&+ r_{N,\,q,\,1},
\end{eqnarray*}
where we sepatated the summand with $m=2q$ above and we used the notation (\ref{rn1}). Consequently,
\begin{eqnarray*}
\mathbf{E} \left[ G_{N} ^ {(2q)}(k, \alpha) ^ {2} \right]&=&\frac{ (c_{2q} ^{k}) ^{2}}{(2q)! }\sum_{v  \in \mathbb{Z}} \varphi_{H, \alpha} (v) ^{2q} -  \frac{ (c_{2q} ^{k}) ^{2}}{(2q)! }\sum_{\vert v\vert \geq N-l+1} \varphi_{H, \alpha} (v) ^{2q} \\
&&+\frac{ (c_{2q} ^{k}) ^{2}}{(2q)! }\sum_{\vert v\vert \leq N-l}\varphi_{H, \alpha} (v) ^{2q} \left(  \frac{N-\vert v\vert -l}{N-l}-1\right)+r_{N,\,q,1}\\
&=& \frac{ (c_{2q} ^{k}) ^{2}}{(2q)! }\sum_{v  \in \mathbb{Z}} \varphi_{H, \alpha} (v) ^{2q} + r_{N,\,q,3}+r_{N,\,q,2}+ r_{N,\,q,1}
\end{eqnarray*}
with
\[r_{N,\,q,2}= \frac{ (c_{2q} ^{k}) ^{2}}{(2q)! }\sum_{\vert v\vert \leq N-l}\varphi_{H, \alpha} (v) ^{2q} \left(  \frac{N-\vert v\vert -l}{N-l}-1\right ),\]
\[r_{N,\,q,3}=  - \frac{ (c_{2q} ^{k}) ^{2}}{(2q)! }\sum_{\vert v\vert \geq N-l+1} \varphi_{H, \alpha} (v) ^{2q}.\]

The asymptotics for $r_{N,\,q,\,1}$ has been studied in Lemma \ref{5d-4}:
\begin{equation*}
\vert r_{N,\,q,\, 1}\vert \leq  C
     \begin{cases}
        \frac{1}{N} & \mbox{ if } p\geq 2 \mbox{ or } p=1, q\geq 2,\\
N ^ {4H-3}  & \mbox{ if } p=1, q=1.
     \end{cases}\\
\end{equation*}
For $r_{N,\,q,2}$ we calculate
\[|r_{N,\,q,\,2}|\leq C\left| \sum_{\vert v\vert \leq N-l}\varphi_{H, \alpha} (v) ^{2q} \left(  \frac{\vert v\vert }{N-l}\right )\right|\leq C \frac{1}{N} \sum_{\vert v\vert \leq N-l} |v|^{2q(2H-2p)+1}.\]
Note that the above series is convergent for $p\geq 2$ or  for  $p=1$ and $q\geq 2$ if $H<1-\frac{1}{2q}$ (which is satisfied for $H<\frac{3}{4}$).  In these cases, we will find the estimate 
\begin{equation*}
\vert r_{N, q, 2}\vert \leq C \frac{1}{N} .
\end{equation*}
For $p=1$ and $q=1$, the sequence $ \sum _{1\leq  v \leq N-l} \vert v\vert ^ {(2H-2p) 2q +1} = \sum _{1\leq \vert v\vert \leq N-l} \vert v\vert ^ { 4H-3}$ behaves as $N ^ {4H-2}$ and we get
\begin{equation*}
\vert r_{N, q, 2}\vert \leq C N ^ {4H-3},
\end{equation*} 
so here we obtain the bounds
\begin{equation*}
\vert r_{N, \,q, 2}\vert \leq  C
     \begin{cases}
        \frac{1}{N} & \mbox{ if } p\geq 2 \mbox{ or } p=1, q\geq 2, H<1-\frac{1}{2q}\\
N ^ {4H-3}  & \mbox{ if } p=1, q=1, H<\frac{3}{4}.
     \end{cases}\\
\end{equation*}
Finally, for $r_{N,\,q,3}$ the same bounds can be established. An application of Lemma \ref{lemma} yields
\begin{equation*}
\vert r_{N, q, 3}\vert \leq C \sum _{\vert v \vert \geq N-l}  \varphi_{H, \alpha} (v) ^{2q}\leq C N ^ {(2H-2p)2q+1} 
\end{equation*}
and consequently,
\begin{equation*}
\vert r_{N,\,q, 3}\vert \leq  C
     \begin{cases}
        \frac{1}{N} & \mbox{ if } p\geq 2 \mbox{ or } p=1, q\geq 2, H<1-\frac{1}{2q}\\
N ^ {4H-3}  & \mbox{ if } p=1, q=1, H<\frac{3}{4}.
     \end{cases}\\
\end{equation*}
Since $\frac{1}{N}<N^{4H-3}$ for $H$ between $\frac{1}{2}$ and $\frac{3}{4}$, the result for $p=1$, $q\geq 2$ follows.
\qed

The following proposition provides a bound for the random part in (\ref{15d-1}).

\begin{prop}\label{propvar}
Let $ G_{N} $ be given by (\ref{gn}). For $q_1,\, q_2\in \{1,\dots,\frac{k}{2}\}$,  $p\geq 2$ and $H\in \left(\frac{1}{2}, 1\right)$, 
\[\operatorname{Var} (\langle DG_{N} ^ {(2q_{1})}(k, \alpha),\, DG_{N} ^ {(2q_{2})}(k, \alpha )\rangle_{\mathcal H})\leq C\frac{1}{N}.\]
For $p=1$ and $H<3\slash 4$
\begin{equation*}
\operatorname{Var} (\langle DG_{N} ^ {(2q_{1})}(k, \alpha),\, DG_{N} ^ {(2q_{2})}(k, \alpha )\rangle_{\mathcal H})\leq C
    \begin{cases}
       \frac{1}{N} & \text{if $H \in (\frac{1}{2}, \frac{5}{8})$ }, \\
       \frac{\log (N)^{3}}{N} & \text{if $H =\frac{5}{8}$}, \\
        N^{8H-6} & \text{if $H \in (\frac{5}{8},\frac{3}{4})$}.
    \end{cases}
\end{equation*}
\end{prop}
{\bf Proof: } We can explicitly compute the Malliavin  derivatives in the statement. For $q \in \{1,\dots,\frac{k}{2}\}$
\[D_\cdotp G_{N} ^ {(2q)}(k, \alpha) = \frac{1}{\sqrt{N-l}} \sum_{i=l} ^{N} \frac{ c_{2q} ^{k}}{(2q-1)! } I_{2q-1} \left(    \frac{ C_{i, \alpha} ^{\otimes (2q-1)}} {( \pi ^{\alpha, N}_{H} (0) ) ^{q}}  \right)C_{i, \alpha} (\cdotp). \]
Assume without loss of generality $q_1\leq q_2$. We have
\begin{eqnarray*}
&& \langle DG_{N} ^ {(2q_{1})}(k, \alpha),\, DG_{N} ^ {(2q_{2})}(k, \alpha )\rangle_{\mathcal H}\\
&=&\frac{1}{(N-l)( \pi ^{\alpha, N}_{H} (0) ) ^{q_1+q_2} } \frac{ c_{2q_1} ^{k}}{(2q_1-1)! } \frac{ c_{2q_2} ^{k}}{(2q_2-1)! } \sum_{i,\,j=l} ^{N}  I_{2q_1-1} \left(    C_{i, \alpha} ^{\otimes (2q_1-1)}  \right) I_{2q_2-1} \left(    C_{j, \alpha} ^{\otimes (2q_2-1)}  \right)\langle C_{i,\,\alpha},\,C_{j,\,\alpha }\rangle_{\mathcal H}\\
&=& \frac{1}{(N-l)( \pi ^{\alpha, N}_{H} (0) ) ^{q_1+q_2} } \frac{ c_{2q_1} ^{k}}{(2q_1-1)! } \frac{ c_{2q_2} ^{k}}{(2q_2-1)! }\\
&&\times \sum_{i,\,j=l} ^{N} \left(\sum_{r=0}^{2q_1-1}r! {2q_1-1 \choose r}{2q_2-1 \choose r} I_{2q_1+2q_2-2-2r}\left(  C_{i, \alpha}^{\otimes (2q_1-1)} \otimes_r   C_{j, \alpha}^{\otimes (2q_2-1)} \right) \right) \langle C_{i,\,\alpha},\,C_{j,\,\alpha }\rangle_{\mathcal H}\\
\end{eqnarray*}
and $\mathbf{E} [\langle DG_{N} ^ {(2q_{1})}(k, \alpha),\, DG_{N} ^ {(2q_{2})}(k, \alpha )\rangle_{\mathcal H}]$ is the term containing $I_0$. It follows that
\begin{eqnarray*}
&& \langle  DG_{N} ^ {(2q_{1})}(k, \alpha),\, DG_{N} ^ {(2q_{2})}(k, \alpha )\rangle_{\mathcal H}-\mathbf{E}[\langle DG_{N} ^ {(2q_{1})}(k, \alpha),\, DG_{N} ^ {(2q_{2})}(k, \alpha )\rangle_{\mathcal H}]\\
&=& \frac{1}{(N-l)( \pi ^{\alpha, N}_{H} (0) ) ^{q_1+q_2} } \frac{ c_{2q_1} ^{k}}{(2q_1-1)! } \frac{ c_{2q_2} ^{k}}{(2q_2-1)! }\\
&&\times \sum_{i,\,j=l} ^{N} \left(\sum_{r=0}^{2q_1-1-w}r! {2q_1-1 \choose r}{2q_2-1 \choose r} I_{2q_1+2q_2-2-2r}\left(  C_{i, \alpha}^{\otimes (2q_1-1)} \otimes_r   C_{j, \alpha}^{\otimes (2q_2-1)} \right) \right) \langle C_{i,\,\alpha},\,C_{j,\,\alpha }\rangle_{\mathcal H},\\
\end{eqnarray*}
where $w=1$ if $l_1\neq l_2$ and $w=2$ otherwise.

Due to the fact that products of integrals of different orders have zero expectation we obtain
\begin{eqnarray*}
P&:=&\mathbf{E}[( \langle  DG_{N} ^ {(2q_{1})}(k, \alpha),\, DG_{N} ^ {(2q_{2})}(k, \alpha )\rangle_{\mathcal H}-\mathbf{E}[\langle  DG_{N} ^ {(2q_{1})}(k, \alpha),\, DG_{N} ^ {(2q_{2})}(k, \alpha )\rangle_{\mathcal H}]])^2]\\
&\sim _{ N\to \infty}&  \frac{C}{(N-l)^2( \pi ^{\alpha, N}_{H} (0) ) ^{2(q_1+q_2)} }\\
&& \times \sum_{r=0}^{2q_1-w}\mathbf{E} \left[\left(\sum_{i,\,j=l}^{N}r!{2q_1-1 \choose r} {2q_2-1 \choose r} I_{2q_1+2q_2-2-2r}(C_{i,\,\alpha}^{\otimes{(2q_1-1)}}\otimes_r C_{j,\,\alpha}^{\otimes{(2q_2-1)}})  \langle C_{i,\,\alpha},\,C_{j,\,\alpha} \rangle_{\mathcal H} \right)^2\right]\\
&=&  \frac{C}{(N-l)^2( \pi ^{\alpha, N}_{H} (0) ) ^{2(q_1+q_2)} } \sum_{r=0}^{2q_1-w}\sum_{i,j,k,m=l}^{N}  r!^2{2q_1-1 \choose r}^2 {2q_2-1 \choose r}^2 \langle C_{i,\,\alpha},\,C_{j,\,\alpha} \rangle_{\mathcal H} \langle C_{k,\,\alpha},\,C_{m,\,\alpha} \rangle_{\mathcal H} \\
&&\times \langle \widetilde{ C_{i,\,\alpha}^{\otimes{(2q_1-1)}}\otimes_r C_{j,\,\alpha}^{\otimes{(2q_2-1)}}},\, \widetilde{ C_{k,\,\alpha}^{\otimes{(2q_1-1)}}\otimes_r C_{m,\,\alpha}^{\otimes{(2q_2-1)}}} \rangle_{\mathcal H^{\otimes(2q_1+2q_2-2-2r)}}=: \sum_{r=0}^{2q_1-w} P_r.
\end{eqnarray*}
We can compute the contractions involved and get via (\ref{contra})
\begin{eqnarray*}
C_{i,\,\alpha}^{\otimes{(2q_1-1)}}\otimes_r C_{j,\,\alpha}^{\otimes{(2q_2-1)}}=   C_{i,\,\alpha}^{\otimes{(2q_1-r-1)}}\otimes C_{j,\,\alpha}^{\otimes{(2q_2-r-1)}} \left(\langle C_{i,\,\alpha},\, C_{j,\,\alpha}\rangle_{\mathcal H}\right)^r.
\end{eqnarray*}
Consequently, one can write for $r\geq 0$
\begin{eqnarray*}
&& \left|\langle \widetilde{ C_{i,\,\alpha}^{\otimes{(2q_1-1)}}\otimes_r C_{j,\,\alpha}^{\otimes{(2q_2-1)}}},\, \widetilde{ C_{k,\,\alpha}^{\otimes{(2q_1-1)}}\otimes_r C_{m,\,\alpha}^{\otimes{(2q_2-1)}}} \rangle_{\mathcal H^{\otimes(2q_1+2q_2-2-2r)}}\right|\\
&\lesssim&\left|\left(\langle C_{i,\,\alpha},\, C_{j,\,\alpha}\rangle_{\mathcal H}\langle C_{k,\,\alpha},\, C_{m,\,\alpha}\rangle_{\mathcal H}\right)^r\right| \\
&&\times \max_{a=0,\dots, l_1-r-1}\left|\langle C_{i,\,\alpha},\, C_{k,\,\alpha}\rangle_{\mathcal H}^{l_1-r-1-a}\langle C_{j,\,\alpha},\, C_{m,\,\alpha}\rangle_{\mathcal H}^{l_2-r-1-a}\langle C_{i,\,\alpha},\, C_{m,\,\alpha}\rangle_{\mathcal H}^a\langle C_{j,\,\alpha},\, C_{k,\,\alpha}\rangle_{\mathcal H}^a\right|
\end{eqnarray*}
due to symmetrisation: the maximum is taken over all outcomes of different permutations of the first and second component of the inner product, the number $a$ signifying the number of $C_i$ in the first component that are appearing in the same places as $C_m$ in the second component in a given permutation.

In total, we obtain for a fixed $r\in \{0,\, 2q_1-w\}$
\begin{eqnarray*}
|P_r|&\leq& C  \frac{1}{(N-l)^2( \pi ^{\alpha, N}_{H} (0) ) ^{2(q_1+q_2)} } \sum_{i,j,k,m=l}^{N} \left|\langle C_{i,\,\alpha},\, C_{j,\,\alpha}\rangle_{\mathcal H}^{r+1}\langle C_{k,\,\alpha},\, C_{m,\,\alpha}\rangle_{\mathcal H}^{r+1}\right|\\
&& \times \max_{a=0,\dots, l_1-r-1}\left|\langle C_{i,\,\alpha},\, C_{k,\,\alpha}\rangle_{\mathcal H}^{l_1-r-1-a}\langle C_{j,\,\alpha},\, C_{m,\,\alpha}\rangle_{\mathcal H}^{l_2-r-1-a}\langle C_{i,\,\alpha},\, C_{m,\,\alpha}\rangle_{\mathcal H}^a\langle C_{j,\,\alpha},\, C_{k,\,\alpha}\rangle_{\mathcal H}^a\right|\\
&=& C \frac{1}{(N-l)^2} \sum_{i,j,k,m=l}^{N}\left|\rho_H^{\alpha ,\, N} (i-j)^{r+1}\rho_H^{\alpha ,\, N} (k-m)^{r+1}\right|\\
&&\times \max_{a=0,\dots, l_1-r-1}\left|\rho_H^{\alpha ,\, N}(i-k)^{l_1-r-1-a}\rho_H^{\alpha ,\, N}(j-m)^{l_2-r-1-a}\rho_H^{\alpha ,\, N}(i-m)^a\rho_H^{\alpha ,\, N}(j-k)^a\right|
\end{eqnarray*}
with $\rho_H^{\alpha ,\, N}$ defined in (\ref{ro}). Due to boundedness of $\rho_H^{\alpha ,\, N}$ we can without loss of generality reduce the number of factors. In particular,
\begin{eqnarray*}
&& \max_{a=0,\dots, l_1-r-1}\left|\rho_H^{\alpha ,\, N}(i-k)^{l_1-r-1-a}\rho_H^{\alpha ,\, N}(j-m)^{l_2-r-1-a}\rho_H^{\alpha ,\, N}(i-m)^a\rho_H^{\alpha ,\, N}(j-k)^a\right|\\
&\leq & C |\rho_H^{\alpha ,\, N}(i-k)\rho_H^{\alpha ,\, N}(j-m)|,
\end{eqnarray*}
since either the factor $|\rho_H^{\alpha ,\, N}(i-k)\rho_H^{\alpha ,\, N}(j-m)|$ or $|\rho_H^{\alpha ,\, N}(i-m)\rho_H^{\alpha ,\, N}(j-k)|$ is contained in the product and for symmetry reasons there is no need to distinguish between these cases. Using this inequality and bounding the first two factors in the same way we arrive at a bound
\begin{eqnarray*}
|P_r|& \leq &C \frac{1}{(N-l)^2} \sum_{i,j,k,m=l}^{N} |\rho_H^{\alpha ,\, N} (i-j)\rho_H^{\alpha ,\, N} (k-m)\rho_H^{\alpha ,\, N}(i-k)\rho_H^{\alpha ,\, N}(j-m)|\\
&\leq& C  \frac{1}{(N-l)^2}N\left(\sum_{v=1}^N |\rho_H^{\alpha ,\, N} (v)|^{4\slash 3}\right)^3,
\end{eqnarray*}
where the last step follows via Young's inequality as in \cite{KT}. The representation $\rho_H^{\alpha ,\, N}(v)=\varphi_H(v)+b_{N,\, H}(v)$ together with the fact that for $|v|\leq N$ we have $b_{N,\, H}(v)\leq C \varphi_H(v)$ for some constant $C$ allows us to replace $\rho_H^{\alpha ,\, N}$ with $\varphi _ H(v)$ in the last bound, since the powers involved are positive. Finally, by Lemma \ref{lemma}
\[ \sum_{v=1}^N |\varphi_H (v)|^{4\slash 3} =   \begin{cases}
       \mathrm{O}(1) & \text{if $H \in (0,\frac{5}{8})$ }, \\
       \mathrm{O}(\log (N)) & \text{if$H =\frac{5}{8}$}, \\
        \mathrm{O}(N^{\frac{8H}{3}-\frac{5}{3}}) & \text{if $H \in (\frac{5}{8},1)$} ,
    \end{cases} \]
 and thus the result follows.
\qed

Before stating and proving the main result of this chapter, let us briefly recall the definition of the Wasserstein distance. The Wasserstein distance between the laws of two $\mathbb{R} ^ {d}$-valued random variables $F$ and $G$ is defined as
\begin{equation}
\label{dw}
d_{W} (F, G)= \sup_{h\in \mathcal{A}}\left| \mathbf{E}h(F)-\mathbf{E}h(G)\right| 
\end{equation}
where $\mathcal{A}$ is the class of Lipschitz continuous function $h:\mathbb{R} ^ {d} \to \mathbb{R}$ such that $\Vert h\Vert _{Lip} \leq 1$, where
$$\Vert h\Vert _{Lip}= \sup_{x, y\in \mathbb{R} ^ {d}, x\not=y} \frac{ \vert h(x)-h(y)\vert} {\Vert x-y\Vert _{\mathbb{R} ^ {d}}}.$$

\begin{theorem}\label{tt2}
Let $p\geq 2$ and consider the sequence (\ref{fnka2}).  Let $Z\sim N(0,1).$ Then there exists a constant $C$ such that
$$d_{W} ( G_{N} (k, \alpha), Z)\leq C \frac{1}{\sqrt{N}}.$$
For $p=1$ and $H<3\slash 4$
\begin{equation*}
d_{W} ( G_{N} (k, \alpha), Z)\leq C
    \begin{cases}
       \frac{1}{\sqrt{N}} & \text{if $H \in (\frac{1}{2}, \frac{5}{8})$ }, \\
       \frac{\log (N)^{3\slash 2}}{\sqrt{N}} & \text{if $H =\frac{5}{8}$}, \\
        N^{4H-3} & \text{if $H \in (\frac{5}{8},\frac{3}{4})$}.
    \end{cases}
\end{equation*}
\end{theorem}
{\bf Proof: } Consider the function $f: \mathbb{R} ^ {\frac{k}{2}}\to \mathbb{R}$, $f(x)= \frac{2}{k} (x_{1}+...+x_{\frac{k}{2}}). $. Note that $f$ is a Lipschitz continuous function with $\Vert f\Vert \leq 1$. From Lemma \ref{det} and Proposition \ref{propvar} it is easy to see that 
$$ d_{W}( \frac{k}{2} {\bf G}_{N} (k, \alpha), \frac{k}{2} {\bf Z} ) =d_{W}\left( \frac{k}{2} (G_{N} ^ {(2)}(k, \alpha),..., G_{N} ^ {(k)}(k, \alpha)), \frac{k}{2} {\bf Z}\right)\leq C \frac{1}{\sqrt{N}}$$
where ${\bf Z}\sim N(0, C)$ if $p\geq 2$ and
$$ d_{W}( \frac{k}{2} {\bf G}_{N} (k, \alpha), \frac{k}{2} {\bf Z} ) \leq C
    \begin{cases}
       \frac{1}{\sqrt{N}} & \text{if $H \in (\frac{1}{2}, \frac{5}{8})$ }, \\
       \frac{\log (N)^{3\slash 2}}{\sqrt{N}} & \text{if $H =\frac{5}{8}$}, \\
        N^{4H-3} & \text{if $H \in (\frac{5}{8},\frac{3}{4})$}.
    \end{cases}$$
if $p=1$ and $H<3\slash 4$.  Now,
\begin{eqnarray*}
d_{W} (G_{N}(k, \alpha) , Z) &=& \sup_{\Vert g\Vert _{Lip} \leq 1} \left| \mathbf{E} g(G_{N} (k, \alpha) )-\mathbf{E}g(Z) \right| \\
&=& \sup_{\Vert g\Vert _{Lip} \leq 1}\left| \mathbf E (g\circ f) (\frac{k}{2}{\bf G} _{N}(k,\,\alpha) ) - \mathbf E (g\circ f) (\frac{k}{2} Z)\right| \\
& \leq & \sup_{\Vert h\Vert _{Lip} \leq 1}\left| \mathbf{E} h(\frac{k}{2}{\bf G}_{N} (k,\,\alpha) )- \mathbf{E}h({\bf Z})\right| \\
&=& d_{W} (\frac{k}{2}{\bf G}_{N}(k,\,\alpha), {\bf Z}).
\end{eqnarray*}\qed

\begin{remark}
For $p=1$ and $k=2$ we retrieve the bounds obtained in \cite{KT}, which also coincide with the speed of convergence for the quadratic variations of the fBm (see \cite{N-P}) under the Wasserstein distance.  For $k=2$, it is also possible to get optimal rates under the total variation distance based on the criteria in \cite{NP3}.
\end{remark}

\subsection{Multivariate Central Limit Theorem}
In this part, we restrict ourselves to the case of quadratic variations (i. e. $k=2$) and we derive a multidimensional CLT. This result will be needed in Section \ref{estimation} which deals with the estimation of the Hurst parameter of the solution to (\ref{systeme wave}).

To establish the multidimensional convergence, we will use Theorem 6.2.3 in \cite{P-T}.
\begin{theorem}\label{multivar}
Let $d \geqslant 2 $ and  $q_{d}, \ldots, q_{1} \geqslant 1$ be some fixed integers. 
Consider vectors
\begin{equation*}
F_{n}=\left(F_{1,n}, \ldots, F_{d,n})=(I_{q_{1}}(f_{1,n}), \ldots, I_{q_{d}}(f_{d,n})\right)
\end{equation*}
with $f_{i,n} \in \mathcal{H}^{\odot q_{i}}$. Let $C$ be a real-valued symmetric non negative definite matrix and let $N \sim \mathcal{N}_{d}(0,C)$. Assume that 
\begin{equation}
\underset{n \to \infty}\lim\mathbf{E}\left( F_{i,n}F_{j,n}\right)=C(i,j) \text{ for } i,\,j\in \{1,\dots ,d\}.
\end{equation}
Then, as $n$ tends to $\infty$, the following two conditions are equivalent:
\begin{itemize}
\item $F_{n}$ converges in law to N,
\item for every $1 \leqslant i \leqslant d$  $F_{i,n}$ converges in law to $\mathcal{N}(0,C(i,i))$.
\end{itemize}
\end{theorem}

We now state and prove the multivariate CLT for the renormalized sequence (\ref{vna}) with $k=2$.
\begin{prop}
Let $P \geqslant 1$ be an integer and $\alpha^{1}, \ldots, \alpha^{P}$ be filters of order $p_1,\dots p_P \geqslant 1$ and lengths $l_1+1,\dots ,l_P+1$ respectively. Let $ V_{N} (2, \alpha)$ be given by (\ref{vna}).
If  $p_1,\dots ,p_P>H+\frac{1}{4}$, we have 
\begin{eqnarray*}
(\sqrt{N} V_{N}(2,\alpha^{1}), \ldots,\sqrt{N} V_{N}(2,\alpha^{P})) \to \mathcal{N} \left(0, \Theta \right),
\end{eqnarray*}
where $(\Theta)_{i,j=1...,P}$ denotes a $P \times P$ matrix with entries given by 
\begin{equation}\label{theta}
\Theta_{n,\,m}= \frac{t^{2}}{8c_{1}(H)^{2}}\sum_{k=l}^{\infty}\left(\sum_{q_{1}=0}^{l_1}\sum_{q_{2}=0}^{l_2}\alpha^{n}_{q_{1}}\alpha^{m}_{q_{2}} \vert k+q_{1}-q_{2}\vert^{2H} \right)^{2} .
\end{equation}
\end{prop}
{\bf Proof: } By (\ref{vnka}) with $k=2$, with $c_{1}(H), c_{2}(H)$ from (\ref{c1}),

\begin{eqnarray*}
\mathbf{E} \left(V_{N}\left(k,\alpha^{n}\right)V_{N}\left(k,\alpha^{m}\right) \right)&=& \frac{N^{4H+2}}{(N-l)^{2}\left(c_{1}(H)N+c_{2}(H)\right)^{2}} \sum_{i,j=l}^{N}\mathbf{E}\left(I_{2}\left( {C_{i,\alpha^{n}}}^{\otimes 2}\right)I_{2}\left(  {C_{j,\alpha^{m}}}^{\otimes 2}\right) \right)\\
&=&\frac{2N^{4H+2}}{(N-l)^{2}\left(c_{1}(H)N+c_{2}(H)\right)^{2}} \sum_{i,j=l}^{N}\langle C_{i,\alpha^{n}}, C_{j,\alpha^{m}}\rangle ^{2}_{\mathcal{H}}.
\end{eqnarray*}
By (\ref{14d-1}), we have for $i,j=l,.., N$
\begin{eqnarray*}
\langle C_{i,\alpha^{n}}, C_{j,\alpha^{m}}\rangle &=& \mathbf{E}\left(U^{\alpha^{n}}\left(\frac{i}{N}\right)U^{\alpha^{m}}\left(\frac{j}{N}\right)\right)\\
&=&\sum_{q_{1}=0}^{l_1}\sum_{q_{2}=0}^{l_2}\alpha^{n}_{q_{1}}\alpha^{m}_{q_{2}}\mathbf{E}\left(u(t,\frac{i-q_{1}}{N})u(t,\frac{j-q_{2}}{N})\right) \\
&=&\sum_{q_{1}=0}^{l_1}\sum_{q_{2}=0}^{l_2}\alpha^{n}_{q_{1}}\alpha^{m}_{q_{2}}\left(\frac{N^{-2H-1}}{2}c_{H} \vert j-i+q_{1}-q_{2} \vert^{2H+1} - \frac{tN^{-2H}}{4} \vert j-i+q_{1}-q_{2} \vert^{2H}\right).\\
\end{eqnarray*}
Plugging this into and using similar computations as in \cite{KT}, we get
\begin{eqnarray*}
&&\mathbf{E} \left(V_{N}\left(k,\alpha^{n}\right)V_{N}\left(k,\alpha^{m}\right) \right)\sim _{N \to \infty}  \frac{2N^{4H+3}}{(N-l)^{2}\left(c_{1}(H)N+c_{2}(H)\right)^{2}}\\
&&\quad\sum_{k=l}^{N}{\left(\frac{N^{-2H-1}}{2}c_{H}\sum_{q_{1}=0}^{l_1}\sum_{q_{2}=0}^{l_2} \alpha^{n}_{q_{1}}\alpha^{m}_{q_{2}}\vert k+q_{1}-q_{2} \vert^{2H+1} - \frac{tN^{-2H}}{4}\sum_{q_{1}=0}^{l_1}\sum_{q_{2}=0}^{l_2}\alpha^{n}_{q_{1}}\alpha^{m}_{q_{2}} \vert k+q_{1}-q_{2} \vert^{2H} \right)}^{2}\\
&&= P_{1} + P_{2} + P_{3}.
\end{eqnarray*}
Using Lemma \ref{lemma}, we get with $p:=\min (p_n,\,p_m)$
\begin{eqnarray*}
P_{1}\sim _{N\to \infty} \frac{c_{1}(H)}{N} \sum_{v=l}^{N} v^{4H-4p},&P_{2}\sim _{N\to \infty}\frac{c_{2}(H)}{N^{2}}\sum_{v=l}^{N} v^{4H-4p+1},& P_{3}\sim _{N \to \infty} \frac{c_{3}(H)}{N^{3}}\sum_{v=l}^{N} v^{4H-4p+2}. \\
\end{eqnarray*}
This shows that $P_{1}$  is the dominant term and it converges for $H<p+\frac{1}{4}$, while the other terms are negligible.  We thus obtain the claimed limit:

\begin{eqnarray*}
\mathbf{E} \left(\sqrt{N}V_{N}\left(k,\alpha^{n}\right)\sqrt{N}V_{N}\left(k,\alpha^{m}\right) \right)&\underset{N \to \infty}\to  & \Theta_{n,\,m},
\end{eqnarray*}
where $\Theta_{n,\,m}$ are given by (\ref{theta}). The second part of the equivalence in Theorem \ref{multivar} was proved as a particular case of the CLT for higher powers, and thus the statement of the proposition follows. \qed

\section{Noncentral limit theorem}\label{nclt}
The asymptotic normality obtained in the previous section holds for any filter of order $p\geq 2$ or for any filter of order $p=1$ and $H<\frac{3}{4}$. It remains to understand what happens  in the case $p=1$ and $H>\frac{3}{4}$. We consider in this section the filter $\alpha =(1, -1)$ (which has order $p=1$) and we will  show that, after a proper normalization, the quadratic variation associated to this filter converges in distribution to a non-Gaussian limit. Let us start by estimating the mean square of the quadratic variation.

\begin{lemma}\label{l1}
Let $ V_{N} (2, (1, -1)) $ be given by (\ref{vna}). If  $v_N:=\mathbf{E} [V_N(2,\,(1,\,-1))^2]$ and $H>\frac{3}{4}$ we have
$$N ^ {4-4H} v_{N} \to _{N \to \infty} \frac{4K_{0}}{k_{1} ^{2}},$$
where $K_{0} $ is given by (\ref{k0}) and $k_{1}$ appears in (\ref{14d-1}).
\end{lemma}
{\bf Proof: } As in Lemma 2 in \cite{KT}, we have by (\ref{14d-1})
\begin{eqnarray*}
 v_{N}&=&\frac{ 2N
^{4H}} { (k_{1}N+ k_{2} ) ^{2}} \sum_{i,j=0}
^{N-1}{\left[\mathbf{E}\left( \Big(u(t,
\frac{i+1}{N})-u(t,\frac{i}{N}) \Big)\Big( u(t,
\frac{j+1}{N})-u(t,\frac{j}{N})\Big) \right)\right]}^2\nonumber\\
&=& \frac{2 N ^{4H}} { (k_{1}N+ k_{2} ) ^{2}} \sum_{i,j=1} ^{N} \left[ k_{1} \frac{\Phi_{H} (i-j) }{ N ^{2H}} + k_{2}\frac{ \Phi_{H+\frac{1}{2}} (i-j)}{ N ^{2H+1} }  \right]^{2}\nonumber \\
&=&  \frac{4 N ^{4H}} { (k_{1}N+ k_{2} ) ^{2}} \sum_{j=1} ^ {N} \sum_{i=j+1} ^ {N-1}  \left[ k_{1} \frac{\Phi_{H} (i-j) }{ N ^{2H}} + k_{2}\frac{ \Phi_{H+\frac{1}{2}} (i-j)}{ N ^{2H+1} }  \right]^{2}\\
&&+  \frac{2 N ^{4H}} { (k_{1}N+ k_{2} ) ^{2}}\sum_{i=1} ^ {N} \left[ k_{1} \frac{1}{N ^ {2H}} + k_{2} \frac{1}{N ^ {2H+1}}\right]^2 \\
&=& \frac{4 N ^{4H}} { (k_{1}N+ k_{2} ) ^{2}} \sum_{l=1} ^{N} \left[ k_{1} \frac{\Phi_{H} (l )}{ N ^{2H}} + k_{2}\frac{ \Phi_{H+\frac{1}{2}} (l)}{ N ^{2H+1} }  \right]^{2}(N-l)\\
&&+  \frac{2 N ^{4H}} { (k_{1}N+ k_{2} ) ^{2}}\sum_{i=1} ^ {N} \left[ k_{1} \frac{1}{N ^ {2H}} + k_{2} \frac{1}{N ^ {2H+1}}\right]^2. 
\end{eqnarray*}
The last summand satisfies
$$ \frac{2 N ^{4H}} { (k_{1}N+ k_{2} ) ^{2}}\sum_{i=1} ^ {N} \left[ k_{1} \frac{1}{N ^ {2H}} + k_{2} \frac{1}{N ^ {2H+1}}\right]^2 \leq C$$
for $N$ large enough while the first summand converges to infinity, see below. Using the asymptotic behavior of $\Phi_{H}$ and $\Phi _{H+\frac{1}{2}}$, namely

$$\Phi_{H}(l)= H(2H-1) l ^ {2H-2} + o (l ^ {2H-2}) $$
and
$$\Phi_{H+\frac{1}{2}} (l)= H(2H+1) l ^ {2H-1}+ o(l ^ {2H-1})$$
for $l$ large, we obtain
\begin{eqnarray*}
 v_{N}&\sim _{N \to \infty}& \frac{4}{k_{1} ^ {2}} N ^ {4H-2}\sum_{l=1} ^ {N}   \left[ k_{1}H(2H-1) \frac{l ^ {2H-2}}{ N ^{2H}} + k_{2}H(2H+1)\frac{  l^ {2H-1}}{ N ^{2H+1} }  \right]^{2}(N-l)\\
&=&  \frac{4}{k_{1} ^ {2}} N ^ {4H-4}
 \frac{1}{N} \sum_{l=1} ^ {N} \left[ k_{1}H(2H-1)\left( \frac{l}{N}\right) ^ {2H-2} + k _{2}H(2H+1)\left( \frac{l}{N}\right) ^ {2H-1} \right] ^ {2}\left(\frac{N-l}{N}\right)
\end{eqnarray*}
and therefore
$$N ^ {4-4H} \frac{ k_{1} ^ {2}}{4K_{0}} v_{N} \to _{N\to \infty} 1$$
with
\begin{eqnarray}K_{0}&=& \int_{0} ^ {1} \left( k_{1} H(2H-1)x ^ {2H-2}+ k_{2}H(2H+1) x ^ {2H-1} \right) ^ {2} (1-x) dx\nonumber \\
&=& k_{1} ^ {2} \frac{ H ^ {2} (2H-1)} {2(4H-3)} + 2k_{1}k_{2} \frac{ H ^ {2} (2H+1)}{2(4H-1)} + k_{2} ^ {2} \frac{ H(2H+1) ^ {2}}{4(4H-1)}.\label{k0}
\end{eqnarray} \qed

Recall that the solution to the wave equation with fractional-white noise can be written as

\begin{equation}\label{sol}
u(t,x)=\int_{0}^{t}\int_{\mathbb{R}^d}G_1(t-s,x-y)W^H(\mathrm{d}s,\mathrm{d}y).
\end{equation}
Let  $x_{i}= \frac{i}{N}$, $i=0,1,.., N$ be a partition of the unit interval $[0,1]$. Denote

\begin{equation*}
g_{t,i} (s,\,x)= G_{1} (t-s, x_{i+1}-x) - G_{1} (t-s, x_{i}-x) 
\end{equation*}
for $i=0,1,.., N-1$ and for $t\geq 0, x\in \mathbb{R}$, with $G_{1}$ given by (\ref{g1}). We can write 

$$u(t, x_{i+1}) - u(t, x_{i}) =  I^{W}_{1} ( g_{t, i } ) $$
where $I^{W} _{1}$ represents the multiple integral of order 1 with respect to the fractional-white Gaussian noise $W ^ {H}$. Then we have

\begin{equation*}
V_{N} (2, (1,-1)):=V_{N}= \frac{1}{N} \sum_{i=1} ^ {N} \frac{ I^{W}_{2} (g_{t,i} ^ {\otimes 2} ) }{ \mathbf{E} \left( u(t, x_{i+1}) - u(t, x_{i})\right) ^ {2}}
\end{equation*}
and so
\begin{equation}
\label{fnn}
 F_{N}= \frac{V_{N}}{\sqrt{v_{N}}}= I_{2} (f_{N}) \mbox{ with } f_{N}(x_{1}, x_{2})= \frac{1}{\sqrt{Nv_{N}}} \frac{N ^ {2H+\frac{1}{2}}}{k_{1}N+k_{2}} \sum_{i=1}^ {N} g_{t,i} ^ {\otimes 2} (x_{1}, x_{2}).
\end{equation}
Since in this part we will use the multiple stochastic integrals with respect to the Gaussian noise $W^{H}$ with covariance (\ref{covW}), let us recall some facts about them. Designate by $\xi$ the set of linear combinations of the simple
functions $\mathds{1}_{\{[0,t]\times A\}},\;t\in [0,T],\;A\in
\mathcal{B}_b(\mathbb{R}^{d})$, the canonical Hilbert space
$\mathcal{H}_{W}$ associated to the field $W^H$, when $H
>\frac{1}{2}$, is defined as the closure of the linear space
generated by $\xi$ with respect to the inner product ${\langle
.,.\rangle}_{\mathcal{H}_{W}}$ which is expressed by:
\begin{equation*}\label{inner-pro-noise}
{\langle \mathds{1}_{\{[0,t]\times A\}},\mathds{1}_{\{[0,s]\times
B\}}\rangle}_{\mathcal{H}_{W}}:=\mathbf{E}(W^H_t(A)W^H_s(B))=\alpha_H\lambda(A
\cap B)\int_{0}^t\int_{0}^s{\mid
u-v\mid}^{2H-2}\mathrm{d}u\mathrm{d}v .
\end{equation*}

The scalar product  in $\mathcal{H}_{W}$ is given by\begin{equation*}\label{inner-pro-function} {\langle
f,g\rangle}_{\mathcal{H}_{W}}=\mathbf{E}(W^H(f)W^H(g))=\alpha_H\int_{0}^T\int_{0}^T\int_{\mathbb{R}^d}f(u,x)g(v,x){\mid
u-v\mid}^{2H-2}\mathrm{d}x\mathrm{d}u \mathrm{d}v.
\end{equation*}
 for every $f,g \in\mathcal{H}_{W}$ such that $\int_{0}^T\int_{0}^T\int_{\mathbb{R}^d}\vert f(u,x)g(v,x)\vert {\mid
u-v\mid}^{2H-2}\mathrm{d}x\mathrm{d}u \mathrm{d}v<\infty$. 

It is possible to represent the Wiener integral with respect to $W_{H}$ as an integral with respect to white noise field with Space-time white noise W via a transfer formula given by (see  \cite{T} for details)  

\begin{equation}\label{trans-form}
\int_{0}^{T} \int_{\mathbb{R}} f(s,y)dW^{H}(s,y)=  \int_{\mathbb{R}}  \int_{\mathbb{R}} \left( \int_{\mathbb{R}} \mathds{1}_{\{[0,t]\}}(u)f(u,x)(u-s)_{+}^{H-\frac{3}{2}}du\right)dW(s,y).
\end{equation}

We will analyze the asymptotic behavior of the sequence $F_{N}$. Since $F_{N}$ belongs to the second Wiener chaos,  its law is completely determined by its cumulants (or equivalently, by its moments). That is, if $F,G$ are elements of the second Wiener chaos then $F$ and $G$ have the same law if and only if they have the same cumulants. Moreover, the convergence of the cumulants to cumulants of element of the second Wiener chaos  implies the convergence in distribution. Let us denote by $k_{m} (F)$, $m\geq 1$, the  $m$th cumulant of a random variable $F$. It is defined as
\begin{equation*}
k_{m}(F)= (-i) ^{n} \frac{\partial ^{n}}{\partial t^{n} }\ln \mathbf{E} (e ^{itF}) | _{t=0}.
\end{equation*}
We have the following link between the moments and the cumulants of $F$: for every $m\geq 1$, 
\begin{equation}\label{mom-cum}
k_{m}(F)= \sum_{\sigma =(a_{1},.., a_{r})\in \mathcal{P}( \{1,..,n\} )}(-1) ^ {r-1} (r-1)! \mathbf{E} X ^ {\vert a_{1}\vert}\ldots \mathbf{E} X ^ {\vert a_{r}\vert}
\end{equation}
if $F\in L^{m}(\Omega)$, where $\mathcal{P} (b)$ is the set of all partitions of $b$. In particular, for centered random variables $F$, we have $k_{1}(F)=\mathbf{E} F, k_{2}(F)= \mathbf{E} F^{2},  k_{3}(F)= \mathbf{E}F^{3},  k_{4}= \mathbf{E} F^{4}- (\mathbf{E} F^{2}) ^{2}.$
In the particular situation when $G=I_{2}(f)$ its cumulants can be computed as (see e.g. \cite{Nourdin}, Proposition 7.2 or \cite{T})
\begin{equation}
k_{m}(G)=  2^{m-1}(m-1)! \int_{ \mathbb{R}  ^{m} }d { u}_{1}\ldots d{ u}_{m} f({ u}_{1},
{ u}_{2}) f({ u}_{2}, { u}_{3})\ldots f({ u}_{m-1}, { u}_{m}) f({ u}_{m},{  u}_{1})
\label{cum2}
\end{equation}
with $u_1,\dots , u_m$ possibly being vectors.

Based on the above formula (\ref{cum2}), we obtain the limit in distribution of (\ref{fnn}).

\begin{prop}\label{pp3}
Let $F_{N}$ be given by (\ref{fnn}) with $H>\frac{3}{4}$. Then the sequence $(F_{N})_{N\geq 1} $ converges in distribution to a random variable $F$ whose law is determined by the cumulants (\ref{cum4}) and (\ref{cum3}).
\end{prop}
{\bf Proof: }
Note first that by the transfer formula (\ref{trans-form})  $W^H(g_{t,\,i})$ has a representation as $W(\tilde{g}_{t,\,i})$ for some  (explicitly known) function $\tilde{g}_{t,\,i}$, where $W$ is a two-dimensional Gaussian noise. Therefore $k_{1}(F_{N})=0, k_{2}(F_{N})=1$ and  the above formula for cumulants (\ref{cum2}) can be applied and we obtain for $m\geq 3$

\begin{eqnarray*}
k_{m}(F_{N})&=&2 ^ {m-1} (m-1) ! \left(  \frac{1}{\sqrt{Nv_{N}}} \frac{N ^ {2H+\frac{1}{2}}}{k_{1}N+k_{2}} \right) ^{ m}\\
&&\int_{\mathbb{R} ^ {m}}  \left( \sum_{j_{1}=1}^ {N} \tilde{g}_{t,j_{1}} ^ {\otimes 2} (x_{1}, x_{2})\right)  \left( \sum_{j_{2}=1}^ {N} \tilde{g}_{t,j_{2}} ^ {\otimes 2} (x_{2}, x_{3})\right) .... \left( \sum_{j_{m}=1}^ {N}\tilde{ g}_{t,j_{m}} ^ {\otimes 2} (x_{m}, x_{1})\right) dx_{1}...dx_{m} \\
&=&2 ^ {m-1} (m-1) ! \left(  \frac{1}{\sqrt{Nv_{N}}} \frac{N ^ {2H+\frac{1}{2}}}{k_{1}N+k_{2}} \right) ^{ m}\\
&&\sum_{j_{1},.., j_{m} =1} ^ {N} \left( \int_{\mathbb{R}} \tilde{g}_{t, j_{1}}(x) \tilde{g}_{t, j_{2}} (x) dx \right)  \left( \int_{\mathbb{R}} \tilde{g}_{t, j_{2}}(x) \tilde{g}_{t, j_{3}} (x) dx \right) .... \left( \int_{\mathbb{R}} \tilde{g}_{t, j_{m}}(x) \tilde{g}_{t, j_{1}} (x) dx \right).
\end{eqnarray*}

We use the {isometry formula for multiple integrals with respect to $W$ as well as the transfer formula in order to get}
\begin{eqnarray*}
\int_{\mathbb{R}} \tilde{g}_{t, j_{1}}(x) \tilde{g}_{t, j_{2}} (x) dx& =&\mathbf{E} \left( u(t, x_{i+1}) - u(t, x_{i}) \right) \left(u(t, x_{j+1}) - u(t, x_{j}) \right)\\
&=& k_{1} \Phi_{H} \left( \frac{i-j}{N} \right) +k_{2}\Phi_{H+\frac{1}{2}}  \left( \frac{i-j}{N} \right), 
\end{eqnarray*}
\begin{equation}
\label{phi}
\Phi_{H}(k)= \frac{1}{2} \left( {\vert k+1 \vert}
^{2H}- 2{\vert k\vert} ^{2H} + \vert k-1\vert
^{2H}\right), \hskip0.3cm k\in \mathbb{R},
\end{equation}
and we obtain
\begin{eqnarray*}
&&k_{m}(F_{N})\\
&=&2 ^ {m-1} (m-1) ! \left(  \frac{1}{\sqrt{Nv_{N}}} \frac{N ^ {2H+\frac{1}{2}}}{k_{1}N+k_{2}} \right) ^{ m}\\
&&\sum_{j_{1},.., j_{m}=1} ^{N} \left[ k_{1} \Phi_{H} \left( \frac{j_{1}-j_{2}}{N} \right) +k_{2}\Phi_{H+\frac{1}{2}}  \left( \frac{j_{1}-j_{2}}{N} \right)\right] .....\left[ k_{1} \Phi_{H} \left( \frac{j_{m}-j_{1}}{N} \right) +k_{2}\Phi_{H+\frac{1}{2}}  \left( \frac{j_{m}-j_{1}}{N} \right)\right].
\end{eqnarray*}

By Lemma \ref{l1}
\begin{eqnarray*}
&&k_{m}(F_{N})\\
&\sim _{N \to \infty}&2 ^ {m-1} (m-1) ! (4K_{0}) ^ {-\frac{m}{2}}N ^ {m} \\
&&\sum_{j_{1},.., j_{m}=1} ^{N} \left[ k_{1} \Phi_{H} \left( \frac{j_{1}-j_{2}}{N} \right) +k_{2}\Phi_{H+\frac{1}{2}}  \left( \frac{j_{1}-j_{2}}{N} \right)\right] .....\left[ k_{1} \Phi_{H} \left( \frac{j_{m}-j_{1}}{N} \right) +k_{2}\Phi_{H+\frac{1}{2}}  \left( \frac{j_{m}-j_{1}}{N} \right)\right]
\end{eqnarray*}
By writing
\begin{equation*}
\Phi_{H} \left( \frac{i-j}{N} \right) = H(2H-1) \int_{\frac{i}{N}}^{\frac{i+1}{N}}  \int_{\frac{j}{N}}^{\frac{j+1}{N}} \vert u-v\vert ^{2H-2} dudv
\end{equation*}
and similarly
\begin{equation*}
\Phi_{H+\frac{1}{2}}  \left( \frac{i-j}{N} \right) = H(2H+1) \int_{\frac{i}{N}}^{\frac{i+1}{N}}  \int_{\frac{j}{N}}^{\frac{j+1}{N}} \vert u-v\vert ^{2H-1} dudv
\end{equation*}
we get, for any $m\geq 3$,
\begin{eqnarray*}
&&k_{m}(F_{N})\\
&\sim _{N \to \infty}&2 ^ {m-1} (m-1) ! (4K_{0}) ^ {-\frac{m}{2}}N ^ {m}\sum_{j_{1},.., j_{m}=1} ^{N} \\
&&\int_{0} ^ {1} \int_{0} ^ {1} dudv \left[ k_{1} H(2H-1)  N ^ {-2H} \vert u-v+j_{1}-j_{2} \vert ^ {2H-2} + k_{2} H(2H+1) N ^ {-2H-1} \vert u-v+j_{1}-j_{2} \vert ^ {2H-1}  \right]\\
&&\ldots \\
&&\ldots \\
&&\int_{0} ^ {1} \int_{0} ^ {1} dudv \left[ k_{1} H(2H-1)  N ^ {-2H} \vert u-v+j_{m}-j_{1} \vert ^ {2H-2} + k_{2} H(2H+1) N ^ {-2H-1} \vert u-v+j_{m}-j_{1} \vert ^ {2H-1}  \right]
\end{eqnarray*}

Next, we write 

$$N ^ {-2H} \vert u-v+j_{m}-j_{1} \vert ^ {2H-2} = N ^ {-2} \left| \frac{j_{1}-j_{2}} {N}\right| ^ {2H-2} \left| 1 + \frac{u-v}{j_{1}-j_{2}} \right| ^ {2H-2}$$
and
$$ N ^ {-2H-1} \vert u-v+j_{1}-j_{2} \vert ^ {2H-1} = N ^ {-2} \left| \frac{j_{1}-j_{2}} {N}\right| ^ {2H-1} \left| 1 + \frac{u-v}{j_{1}-j_{2}} \right| ^ {2H-1}$$
and obtain
\begin{eqnarray*}
k_{m}(F_{N})
&\sim _{N \to \infty} &2 ^ {m-1} (m-1) ! (4K_{0}) ^ {-\frac{m}{2}}N ^ {-m}\sum_{j_{1},.., j_{m}=1} ^{N} \\
&&\int_{0} ^ {1} \int_{0} ^ {1} dudv \left[ k_{1} H(2H-1)  \left| \frac{j_{1}-j_{2}} {N}\right| ^ {2H-2} \left| 1 + \frac{u-v}{j_{1}-j_{2}} \right| ^ {2H-2}\right. \\
&&\left. + k_{2} H(2H+1)  \left| \frac{j_{1}-j_{2}} {N}\right| ^ {2H-1} \left| 1 + \frac{u-v}{j_{1}-j_{2}} \right| ^ {2H-1}\right]\\
&&\ldots \\
&&\ldots \\
&&\int_{0} ^ {1} \int_{0} ^ {1} dudv \left[ k_{1} H(2H-1)  \left| \frac{j_{m}-j_{1}} {N}\right| ^ {2H-2} \left| 1 + \frac{u-v}{j_{m}-j_{1}} \right| ^ {2H-2}\right. \\
&&\left. + k_{2} H(2H+1)  \left| \frac{j_{m}-j_{1}} {N}\right| ^ {2H-1} \left| 1 + \frac{u-v}{j_{m}-j_{1}} \right| ^ {2H-1}\right].
\end{eqnarray*}
We claim that
\begin{eqnarray*}
k_{m}(F_{N})
&\sim _{N \to \infty}&2 ^ {m-1} (m-1) ! (4K_{0}) ^ {-\frac{m}{2}}N ^ {-m}\sum_{j_{1},.., j_{m}=1} ^{N} \\
&& \left[ k_{1} H(2H-1)  \left| \frac{j_{1}-j_{2}} {N}\right| ^ {2H-2}+  k_{2} H(2H+1)  \left| \frac{j_{1}-j_{2}} {N}\right| ^ {2H-1} \right]\\
&&\ldots \\
&&\ldots\\
&& \left[ k_{1} H(2H-1)  \left| \frac{j_{m}-j_{1}} {N}\right| ^ {2H-2}+  k_{2} H(2H+1)  \left| \frac{j_{m}-j_{1}} {N}\right| ^ {2H-1} \right].
\end{eqnarray*}
This follows by a standard procedure (see \cite{Nourdin} or \cite{T}) from the Taylor expansion in the vicinity of $x=0$  of the functions 
$$1- (1+x) ^ {2H-2} \mbox{ and } 1-(1+x) ^ {2H-1}$$
and by the dominated convergence theorem. Therefore, for $m\geq 3$

\begin{eqnarray*}
k_{m}(F_{N}) &\to  _{N\to \infty} &2 ^ {m-1} (m-1) ! (4K_{0}) ^ {-\frac{m}{2}}\int_{[0,1] ^ {m}}dx_{1}...dx_{m}\\
&&\left( k_{1} H(2H-1) \vert x_{1}-x_{2}\vert ^ {2H-2} +  k_{2} H(2H+1) \vert x_{1}-x_{2}\vert ^ {2H-1}\right)\\
&&\ldots \\
&&\left( k_{1} H(2H-1) \vert x_{m}-x_{2}\vert ^ {2H-2} +  k_{2} H(2H+1) \vert x_{m}-x_{1}\vert ^ {2H-1}\right)
\end{eqnarray*}
Notice that the above integral is finite by Lemma 3.3 in \cite{BaiTa}. Also, clearly
$$k_{1}(F_{N})=0 \mbox{ and } k_{2} (F_{N})=1.$$
Since $F_{N}$ belongs to the second Wiener chaos, the convergence of cumulants determines the convergence of $F_{N}$ in law to a random variable $F$  with cumulants, for $m\geq 3$

\begin{eqnarray}
k_{m}(F)&=&2 ^ {m-1} (m-1) ! (4K_{0}) ^ {-\frac{m}{2}}\int_{[0,1] ^ {m}}dx_{1}...dx_{m}\nonumber \\
&&\left( k_{1} H(2H-1) \vert x_{1}-x_{2}\vert ^ {2H-2} +  k_{2} H(2H+1) \vert x_{1}-x_{2}\vert ^ {2H-1}\right)\nonumber \\
&&\ldots \nonumber  \\
&&\left( k_{1} H(2H-1) \vert x_{m}-x_{1}\vert ^ {2H-2} +  k_{2} H(2H+1) \vert x_{m}-x_{1}\vert ^ {2H-1}\right)\label{cum3}
\end{eqnarray}
and
\begin{equation}
\label{cum4}
k_{1}(F)=0 \mbox{ and } k_{2} (F)=1.
\end{equation}
The existence of such a limit is ensured by the Fr\'echet-Shohat theorem: It follows from the convergence of cumulants that also all the moments of $F_N$ converge to some real numbers $M_m$, $m\in\mathbb N$, as $N$ tends to infinity. Moreover, by hypercontractivity (\ref{hyper}) the $m$th absolute moments of $F_N$ are bounded by $(m-1)^m$. Therefore, also the limits of the moment sequences will be bounded by $(m-1)^m$, which means that the growth condition $\limsup _{m\to\infty}(\frac{1}{m!}|M_m|)^{1\slash m}<\infty$ is satisfied, thus yielding the existence of a limiting distribution with the cumulants obtained above. 
\qed

\vskip0.1cm
Notice that the limit law with cumulants (\ref{cum3}) and (\ref{cum4}) is related to the Rosenblatt distribution but is more complex. For instance, if the constant $k_{2}$ will vanish in (\ref{cum3}), then we would obtain a Rosenblatt distribution in the limit.

\section{Estimation of the Hurst parameter $H$ }\label{estimation}

We will apply the theoretical results from Section \ref{clt} in order to construct and analyze several estimators for the Hurst index of the mild solution (\ref{sol-wave-WFrac}) to the wave equation (\ref{systeme wave}).  It is worth to emphasize that the estimators are based on the observations of the process $u$ at only  a fixed time and at discrete points in space.

We will define two kinds of estimators for the Hurst parameter. For  the first kind  we will consider the observation time $t$ of our equation to be known, and  the estimators obtained will be asymptotically normal with the rate  of convergence of order $\sqrt{N}\log (N)$ for $H<p-\frac{1}{4}$. In the second case we develop an estimator for $H$ if the time $t>1$ is not known. This estimator will also be asymptotically normal, but with a slower rate of convergence, namely $\sqrt{N}$. Both kinds of estimators are strongly consistent.

\subsection{Estimators for known $t$}
We follow the standard procedure from \cite{coeur} or \cite{CTV} to define our estimators. Define the $k$-th empirical absolute moment of discrete varations of the mild solution $u(t,x)$ for fixed time  $t >1$ and a filter $\alpha$.  That is, let

\begin{eqnarray}
S_{N}(k,\alpha)&=& \frac{1}{N-l}\sum_{i=l}^{N-1}\left\vert U^{\alpha}\left(\frac{i}{N}\right) \right\vert^{k} 
\end{eqnarray}
with $U^{\alpha}\left(\frac{i}{N}\right)  $ defined in (\ref{uin}).  Since $U^{\alpha}\left(\frac{i}{N}\right)  $ is Gaussian, we have \(\mathbf{E}\left[\left| U^{\alpha}\left(\frac{i}{N}\right)\right|^k\right]=\left(\pi _H^{\alpha,\, N}(0) \right)^{\frac{k}{2}} E_k,\) where $E_k$ denotes the $k$-th absolute moment of the standard Gaussian random variable, and therefore we obtain
\[\mathbf{E} [S_N (k,\,\alpha)]=\left(\pi _H^{\alpha,\, N}(0) \right)^{\frac{k}{2}} E_k.\]
Thus, for a given $k$,  replacing $\mathbf{E} [S_N (k,\,\alpha)]$ by  $S_N (k,\,\alpha)$, we obtain an estimator for $H$ that is a pointwise solution to the equation
\[S_N (k,\,\alpha)^{\frac{2}{k}}-E_k^{\frac{2}{k}}\pi _x^{\alpha,\, N}(0) =0 \]
with respect to $x$. Recall that (see (\ref{14d-1}))
\[\pi _x^{\alpha,\, N}(0)=\frac{t}{2N^{2x}}\Phi_{x,\,\alpha }(0)-\frac{c_x}{N^{2x+1}}\Phi_{x+\frac{1}{2},\,\alpha}(0)\]
and we denote
\[c_1(x):=\Phi_{x,\,\alpha} (0)=-\frac{1}{2}\sum_{q,\,r=0}^l \alpha_q\alpha_r|q-r|^{2x},\qquad c_2(x)=c_x\Phi_{x+\frac{1}{2},\,\alpha}(0) .\]
Note that for large $N$ the function $g(x):= \pi _x^{\alpha,\, N}(0)$ is invertible. In order to see this we consider the derivative
\[g'(x)=\frac{t}{2}\left(\frac{c_1'(x)}{N^{2x}}-\frac{2\log (N)c_1(x)}{N^{2x}}\right)-\left(\frac{c_2'(x)}{N^{2x+1}}-\frac{2\log (N)c_2(x)}{N^{2x+1}}\right).\]
As shown in \cite{coeur}, the first expression in the parenthesis becomes negative for large $N$, and since it is the asymptotically dominating term, also the whole function will become negative for $N$ large enough. Therefore, for such $N$ the function $g$ is strictly decreasing and we can define estimators by inverting it:
\begin{equation}
\label{h1}
\widehat{H}_{k,\,N}:= \left(\pi _\cdotp^{\alpha,\, N}(0)\right)^{-1}\left(\left(\frac{S_N (k,\,\alpha)}{E_k}\right)^{\frac{2}{k}} \right).
\end{equation}

Another estimator can be obtained by inverting only the dominant part of the empirical absolute moment. Notice that asymptotically $\pi_x ^{\alpha,\, N}(0)$ is equal to $\frac{t}{2N^{2x}}\Phi_{x,\,\alpha }(0)=:\bar{g}(x)$, which is easier to invert than its exact counterpart. This motivates the definition of another class of estimators,

\begin{equation}
\label{h2}
\bar{H}_{k,\,N}:= \bar{g}^{-1}\left(\left(\frac{S_N (k,\,\alpha)}{E_k}\right)^{\frac{2}{k}} \right).
\end{equation}
We show that the two estimators constructed above are consistent and we give their limit behavior in distribution.

\begin{prop}
The estimators $ \widehat{H}_{N,\,k}$and  $\bar{H}_{N,\,k}$  given by (\ref{h1}) and (\ref{h2}) of the Hurst parameter $H>\frac{1}{2}$ are strongly consistent.
Moreover, with $v_N^{(k)}:=\mathbf{E} [V_N(k,\,\alpha)^2]$, for  $H\leq p-\frac{1}{4}$ we have 
\begin{equation*}
\frac{k \log(N)}{\sqrt{v_{N}^{(k)}}}\left(H - \widehat{H}_{N,\,k} \right) \stackrel{ d }{\to} \mathcal{N}(0,1)
\end{equation*}
and for $H>\frac{3}{4}$, $\alpha = (1,\,-1)$, $k=2$
\begin{equation*}
\frac{2 \log(N)}{\sqrt{v_{N}^{(2)}}}\left(H - \widehat{H}_{N,\,2} \right) \stackrel{ d }{\to}_{N \to \infty}F
\end{equation*}
where  $F$ is the random variable from Proposition \ref{pp3}.
The same statements hold for $\bar{H}$.
\end{prop}
{\bf Proof: }Since for every $k\geq 2$,
\begin{equation}\label{17d-1}
v_N^k=
     \begin{cases}
       O(1\slash N) & \text{if $H <p-\frac{1}{4}$ }, \\
       O\left(\frac{\log(N)}{N}\right) & \text{if $ H =\frac{3}{4},\, p=1$}, \\
       O\left(\frac{1}{N^{2-2H}}\right) & \text{if $H>\frac{3}{4},\, \alpha=(1,\,-1),\,k=2$},
     \end{cases}\\
\end{equation}
the almost sure convergence to zero of $V_N$ follows by hypercontractivity with a Borel-Cantelli argument, see e.g. \cite{TTV}. Due to the fact that the functions $g$ and $\bar{g}$ are asymptotically equal we obtain the asymptotic equality of $\bar{H}_{N,k}$ and $\widehat{H}_{N,k}$ and thus also strong consistency of $\bar{H}_{N,k}$.

For the asymptotic behaviour we can refer to the calculations from \cite{coeur} and obtain
\[ V_{N}(k, \alpha):=V_N=k\log (N)(H-\widehat{H}_{N,\,k})(1+o(1)),\]
which means that by Slutsky's Lemma we will get
\begin{equation*}
\frac{k \log(N)}{\sqrt{v_{N}^{(k)}}}\left(H - \widehat{H}_{N,\,k} \right) \stackrel{ d }{\to} \mathcal{N}(0,1)
\end{equation*}
for $H\leq p-\frac{1}{4}$. For $H<p-\frac{1}{4}$ this implies in particular that
\begin{equation*}
k \log(N)\sqrt{N}\left(H - \widehat{H}_{N,\,k} \right) \stackrel{ d }{\to} \mathcal{N}(0,\sigma^2)
\end{equation*}
with $\sigma^2$ defined in Lemma \ref{5d-4}. For $H>\frac{3}{4},\, \alpha=(1,\,-1),\,k=2$ the relation yields
\begin{equation*}
\frac{2 \log(N)}{\sqrt{v_{N}}}\left(H - \widehat{H}_{N,\,2} \right) \stackrel{ d }{\to}F
\end{equation*}
for $F$ given above. The same results follow for $\bar{H}$ due to its asymptotic equality to $\widehat{H}_{N,k}$.\qed

\begin{remark}

Note that this result  provides the following speeds of convergence (see \ref{17d-1}):  $\sqrt{N}\log (N)$ for $H<p-\frac{1}{4}$, $\sqrt{N}\sqrt{\log(N)}$ for $ H =\frac{3}{4},\, p=1$ and $N^{2-2H}\log(N)$ for  $H>\frac{3}{4},\, \alpha=(1,\,-1),\,k=2$.
\end{remark}

\subsection{An estimator for unknown $t$}

Assume that the time $t>1$ at which the solution (\ref{sol-wave-WFrac}) is observed is not known. Similarly to \cite{istas}, if two sequences $(a^{(1)}_i)_{i\in \{0,\dots ,p \}}$ and  $(a^{(2)}_i)_{i\in \{0,\dots ,2p \}}$ are considered, where $a^{(2)}$ is obtained by ''thinning'' the sequence $a^{(1)}$ (i.e., $a^{(2)}_{2k}:= a^{(1)}_k$ for $k\in \{0,\dots ,p \} $ and zero otherwise), then it follows that
\[\Phi_{H,\,a^{(2)}}(0)=2^{2H}\Phi_{H,\,a^{(1)}}(0)\text{ and } \Phi_{H+\frac{1}{2},\,a^{(2)}}(0)=2^{2H+1}\Phi_{H+\frac{1}{2},\,a^{(1)}}(0),\]
which implies that for large $N$ we have approximately $\pi _H^{a^{(2)},\, N}(0) \sim 2^{2H}\pi _H^{a^{(1)},\, N}(0)$. This, in turn, can be transferred to $S_N$:
\[\mathbf{E}[S_N (k,\,a^{(2)})]=\left(\pi _H^{a^{(2)},\, N}(0) \right)^{\frac{k}{2}} E_k\sim 2^{Hk}\left(\pi _H^{a^{(1)},\, N}(0) \right)^{\frac{k}{2}} E_k=2^{Hk}\mathbf{E} [S_N (k,\,a^{(1)})].\]
This motivates another estimator for $H$ defined by
\begin{equation}
\label{h3}
\tilde{H}_N:=\frac{1}{k}\log_2 \left(\frac{S_N(k,\,a^{(2)})}{S_N(k,\,a^{(1)})} \right).
\end{equation}

Its limit behavior is given below.

 \begin{prop}
The estimator $\tilde{H}_N$ (\ref{h3}) is strongly consistent for all $H>\frac{1}{2}$. Moreover for $H < p- \frac{1}{4}$, we have 
\[ \sqrt{N}(\tilde{H}_N-H)\stackrel{d}{\to}N(0,\,\sigma^2)\]
with $\sigma>0$.
\end{prop}
{\bf Proof: }It follows from the fact that $V_N\to _{N\to \infty} 0$ almost surely that $S_N$ converges almost surely to its expectation. Thus, strong consistency is clear by construction of $\tilde{H}$. The multivariate convergence statement yields asymptotic normality by the delta method, similarly to \cite{coeur}.\qed

\subsection{Numerical computations and simulation experiments }
In this section we conduct simulations of the soultion process and compare numerical performances of different estimators introduced in the previous section. More specifically, we are going to analyse the behaviour of $\bar{H}_{2,\,N}$ for filters $(1,\,-1)$ as well as $(1,\,-2,\,1)$, that of its exact counterpart $\hat{H}_{2,\,N}$ for the second filter as well as that of $\tilde{H}_N$ for different values of $H$. For $N=1000$ and $t=3$ we get the following results from MSE computed from $100$ iterations:
\begin{center}
    \begin{tabular}{ c | c | c | c |}
     & $H=0.51$ & $H=0.7$ & $H=0.95$  \\ \hline
    $\bar{H}_{2,\,N} (1,\,-1)$ & $1.02\cdotp 10^{-5}$ & $1.61\cdotp 10^{-5}$ & $0.001$  \\ \hline
    $\bar{H}_{2,\,N}(1,\,-2,\,1)$ & $1.2\cdotp 10^{-5}$ & $9.626\cdotp 10^{-6}$ & $1.98\cdotp 10^{-6}$ \\ \hline
    $\hat{H}_{2,\,N}(1,\,-2,\,1)$ & $1.2\cdotp 10^{-5}$ & $9.634\cdotp 10^{-6}$ & $1.99\cdotp 10^{-6}$ \\ \hline
    $\tilde{H}_N$ & $0.002$ & $0.001$ & $0.001$  \\
    \hline
    \end{tabular} .
\end{center}

The estimator $\tilde{H}_N$ performs the worst. This can be explained heuristically by the fact that it contains two sources of error instead of one, this being the practical trade-off in the case where time $t$ is not available. Another interesting observation is that the exact estimator $\tilde{H}$ is not performing better than the estimator $\bar{H}$ which uses the inverse of an approximation of the actual function, which justifies the use of the simpler version in applications.

\begin{center}
\begin{table}[ht]
    \begin{tabular}{ |c | c | c | c | c|}
  \hline
 True value H  & Mean $\bar{H}_{2,\,N} (1,\,-1)$ &  Mean $\bar{H}_{2,\,N}(1,\,-2,\,1)$ &  Mean $\hat{H}_{2,\,N}(1,\,-2,\,1)$ &  Mean$\tilde{H}_N$ \\
  \hline
 0.51 &0.5107118 & 0.5138851 & 0.5110081 & 0.5110082 \\
 0.55 & 0.5499827 & 0.5362797 & 0.549677 & 0.549678\\
 0.60 & 0.5997487 & 0.6007376 & 0.5999698 & 0.5999722\\
 0.65 & 0.6498786 & 0.6510065 & 0.6502865 & 0.6502909 \\
 0.70 & 0.7005558 & 0.6925 & 0.7003125 & 0.7003196 \\
 0.75 &  0.7500486 &  0.7482407 & 0.7499587 & 0.74997 \\
 0.80 & 0.8005769 & 0.7966326 & 0.7998019 & 0.7998186 \\
 0.85 & 0.8512704 & 0.8517664 &  0.8500505 & 0.8500754 \\
 0.90 & 0.9042009 & 0.8927607 & 0.8997257 & 0.8997638 \\
 0.95 & 0.9587621 & 0.9540507 & 0.9498974 & 0.9499602 \\
 0.99 & 1.01826 &  0.9959974 & 0.9898137 & 0.9899168 \\
  \hline
\end{tabular}
\caption{Mean of the estimators for $100$ simulations}
\end{table}
\end{center}

\begin{figure}[H]
\centering
\includegraphics[width=0.3\textwidth]{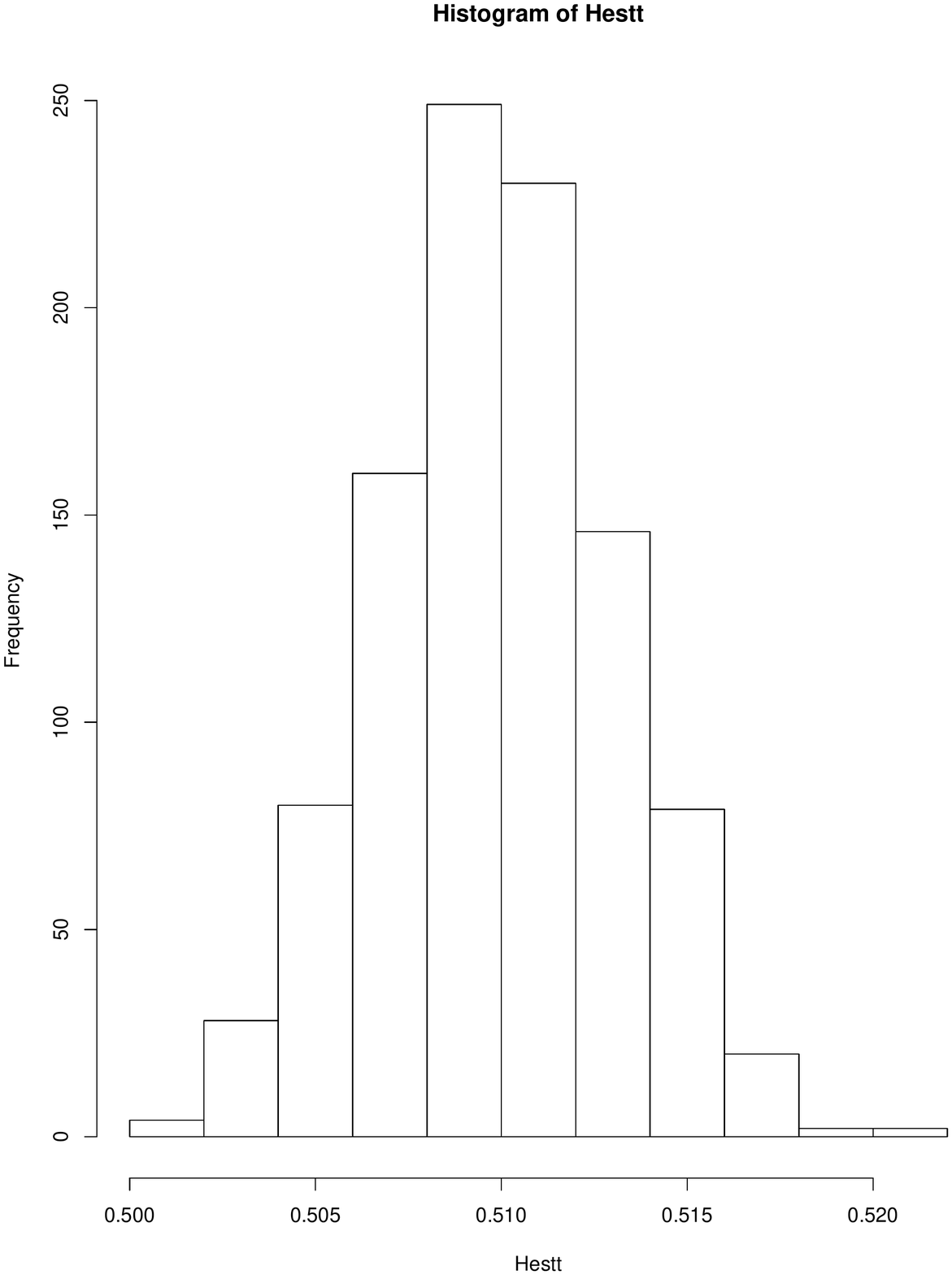}
\includegraphics[width=0.3\textwidth]{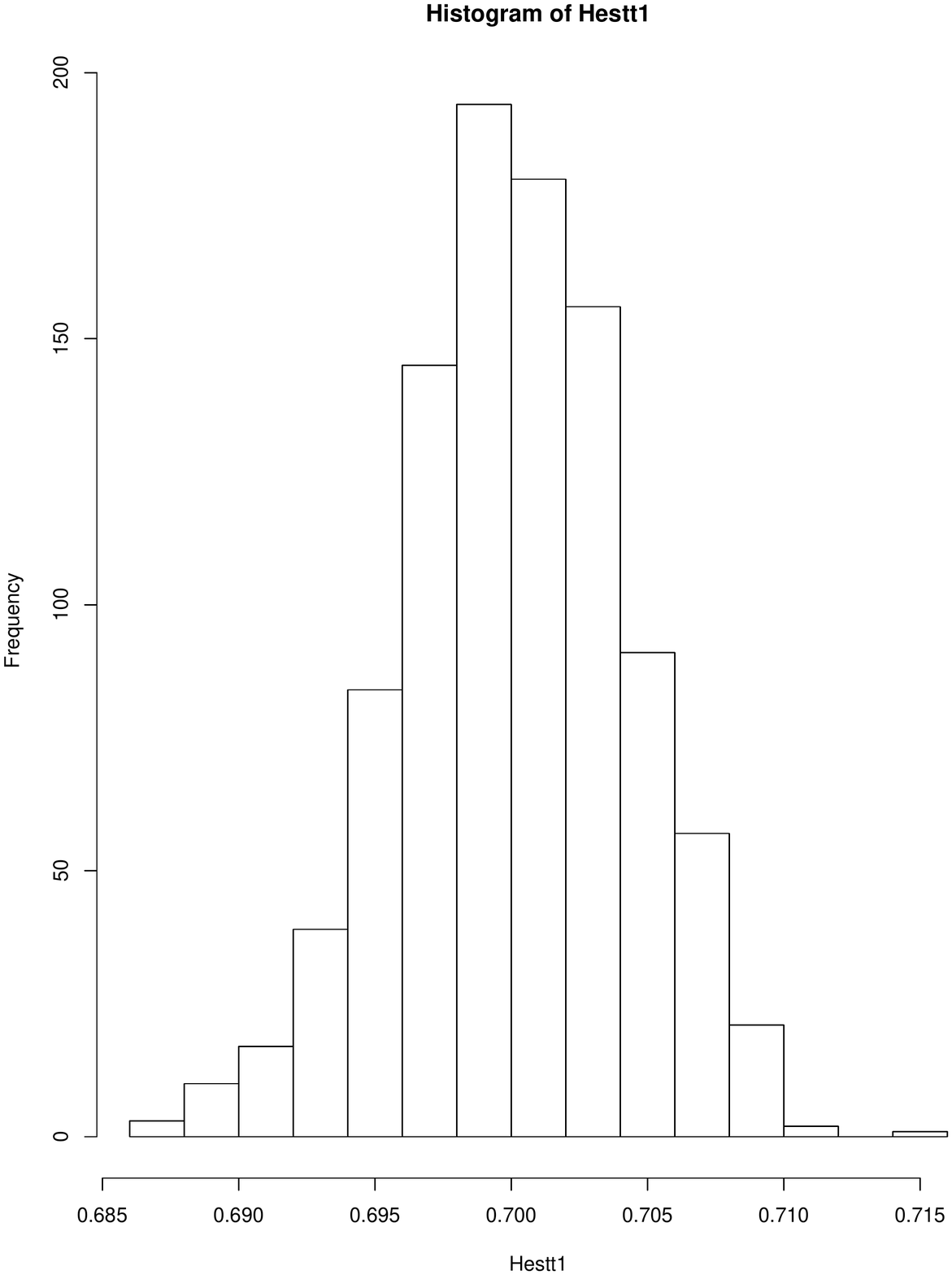}
\includegraphics[width=0.3\textwidth]{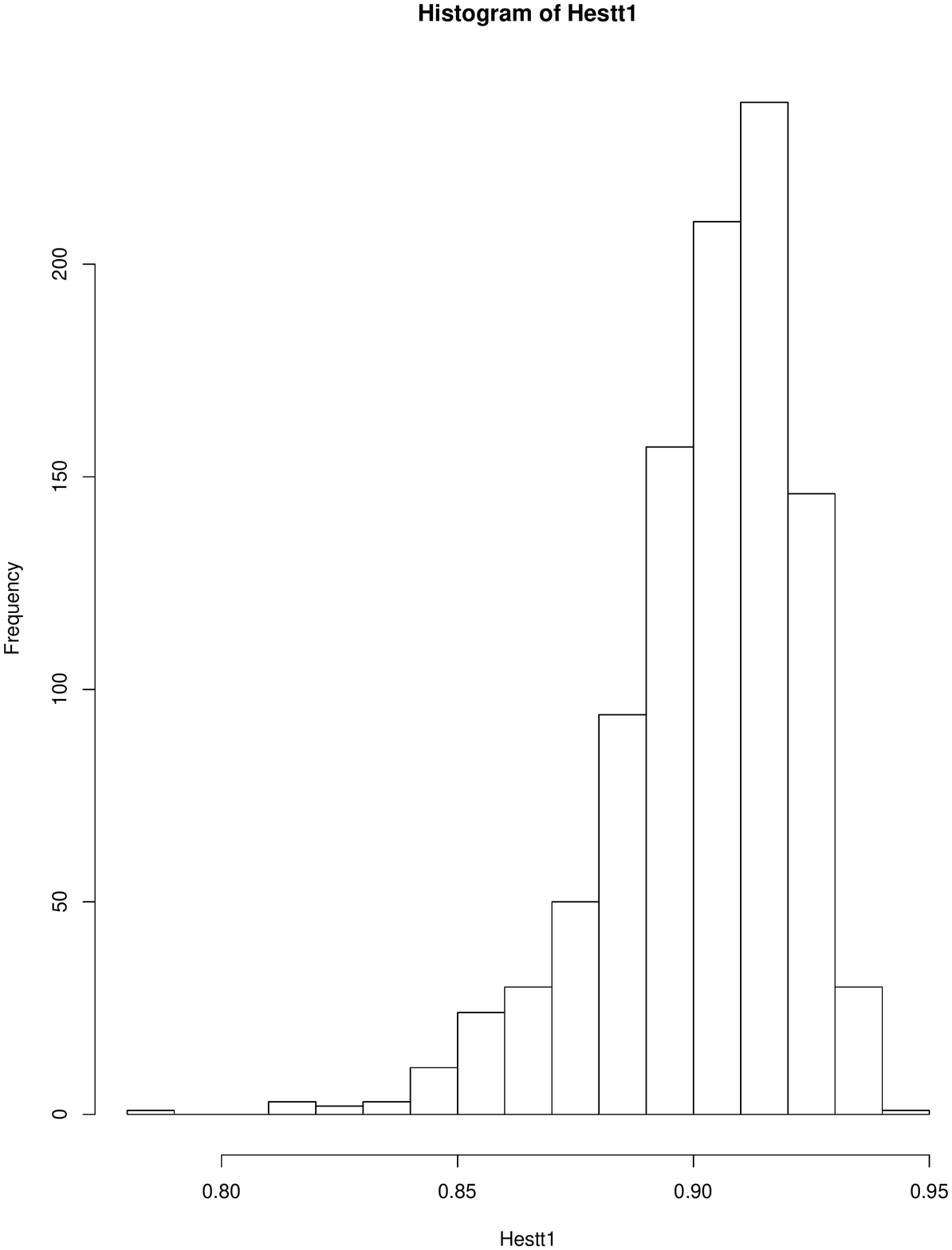}
\caption{Histograms for $H=0.51$, $0.7$ and $0.9$ respectively.}
\label{fig:paths}
\end{figure}
\begin{figure}[H]
\centering
\includegraphics[width=0.3\textwidth]{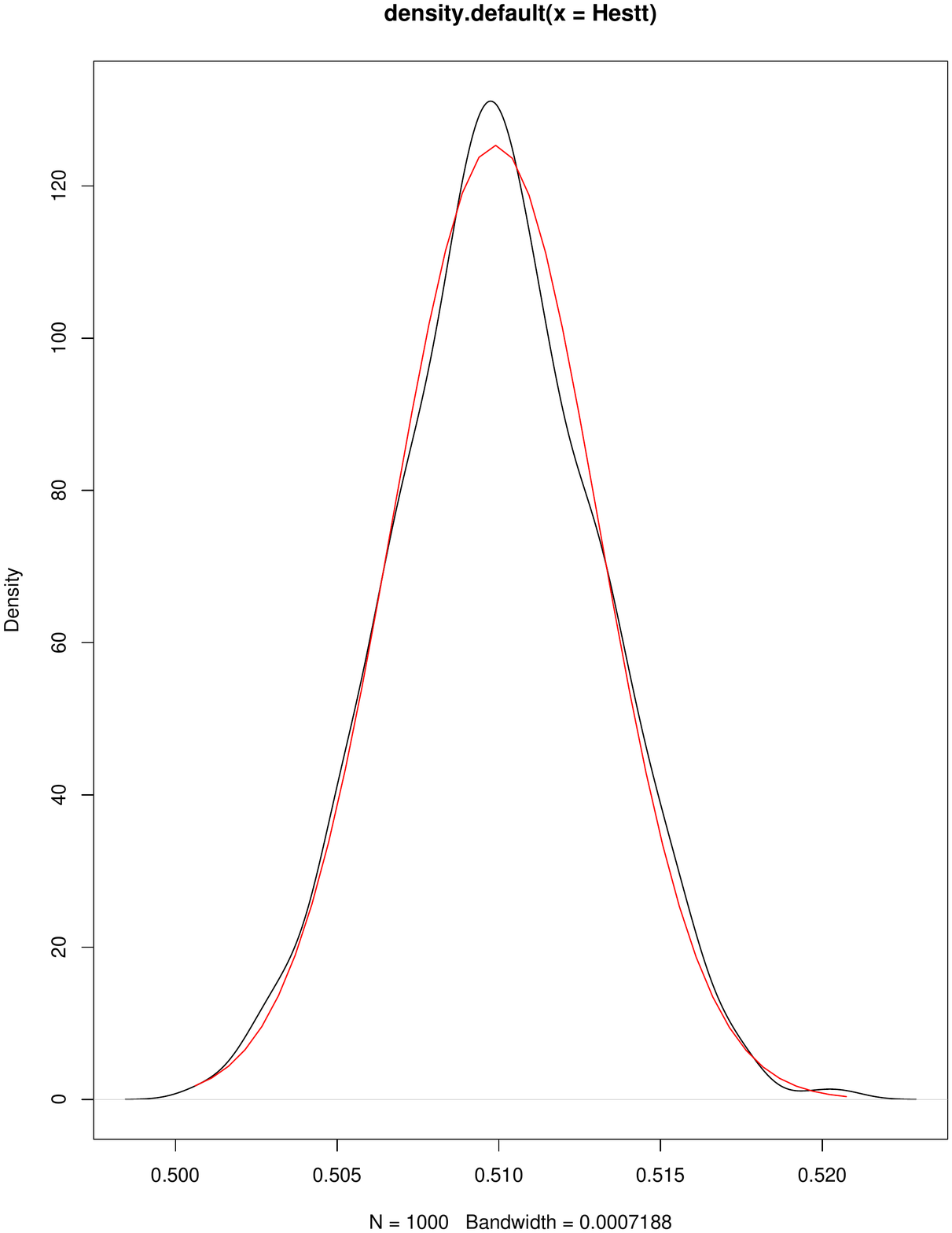}
\includegraphics[width=0.3\textwidth]{densityFit_H051_t3}
\includegraphics[width=0.3\textwidth]{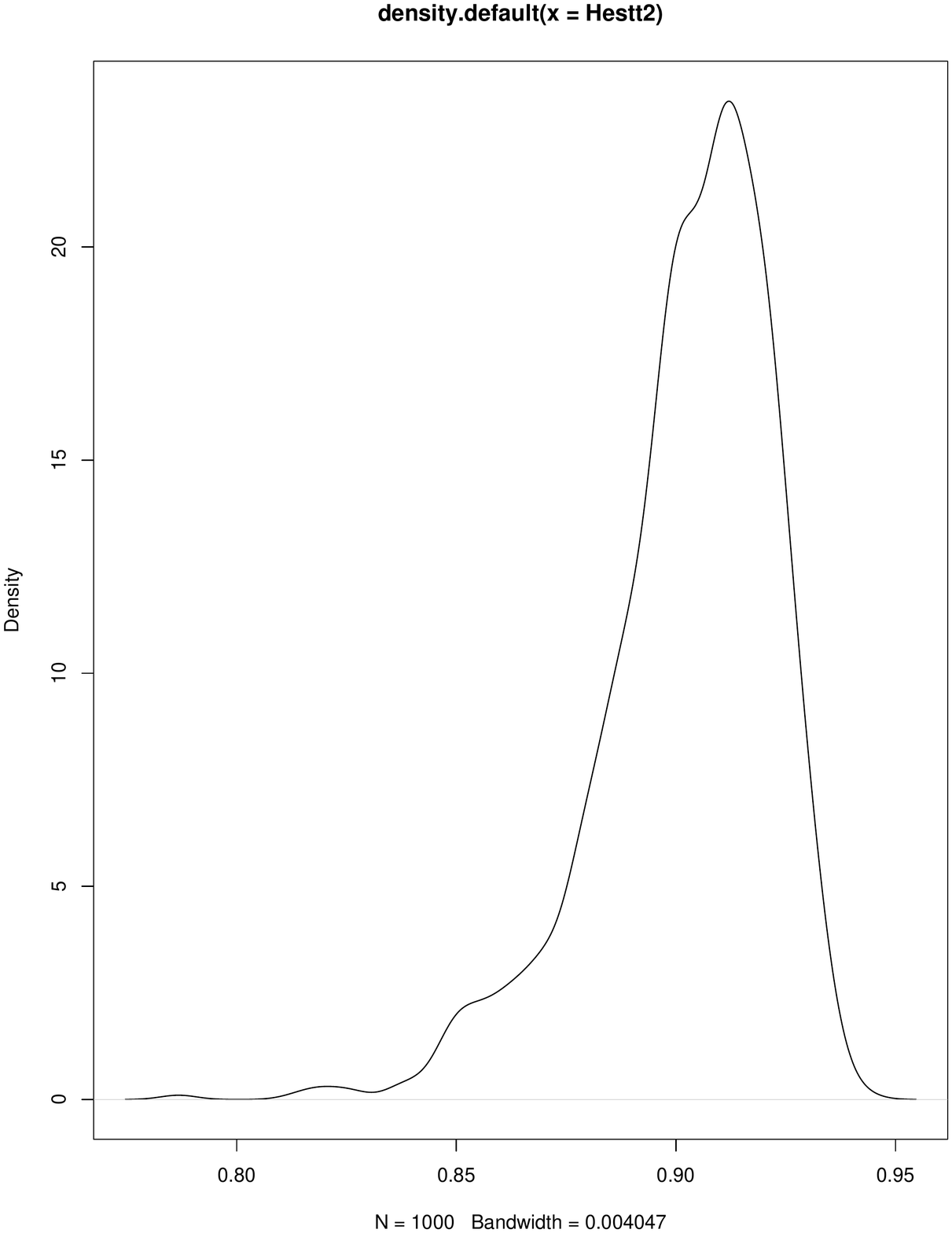}
\caption{Normal fits of empirical densities for $H=0.51$, $0.7$, density plot for $0.9$.}
\label{fig:paths2}
\end{figure}
\begin{figure}[H]
\includegraphics[width=.22\linewidth]{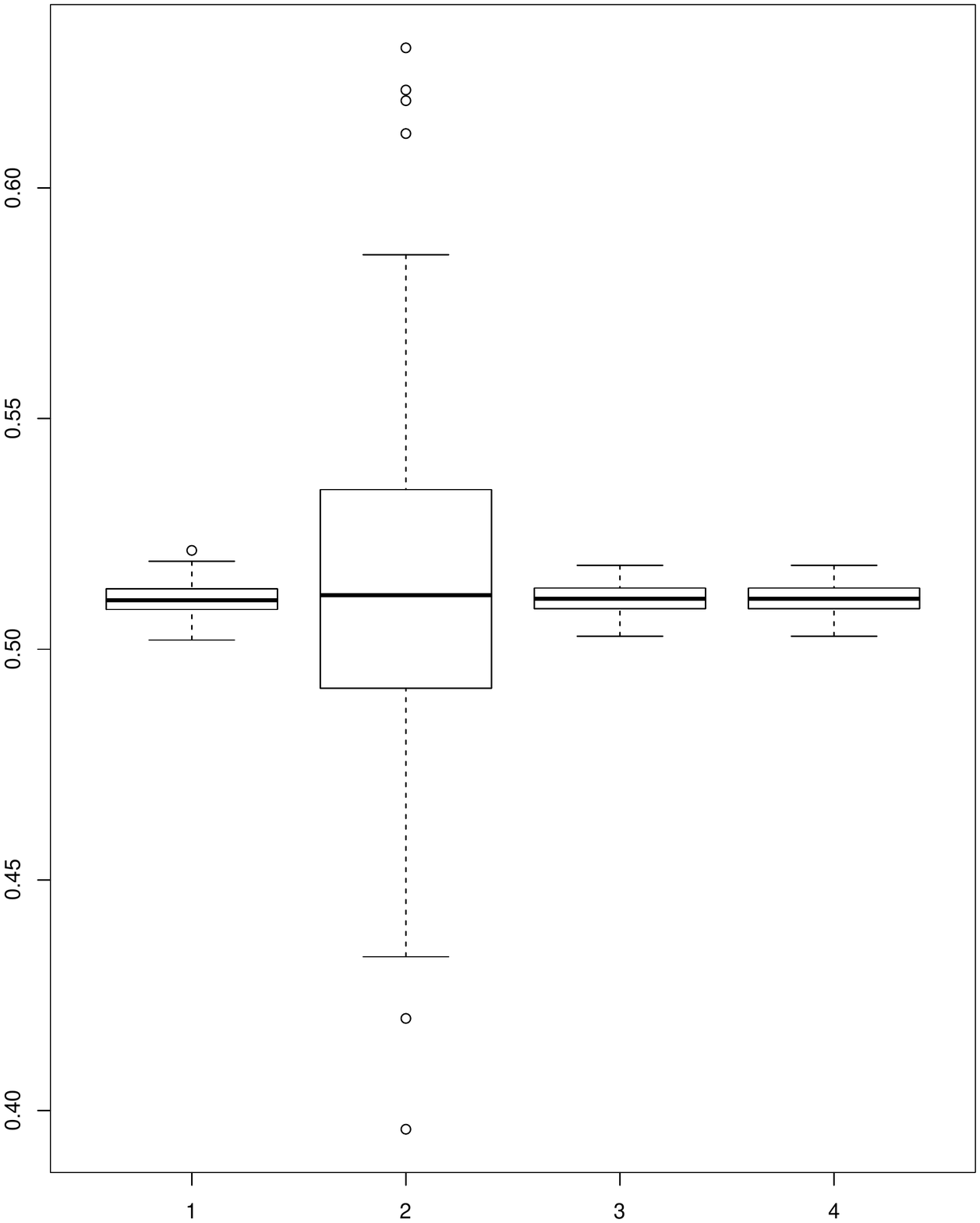}\quad\includegraphics[width=.22\linewidth]{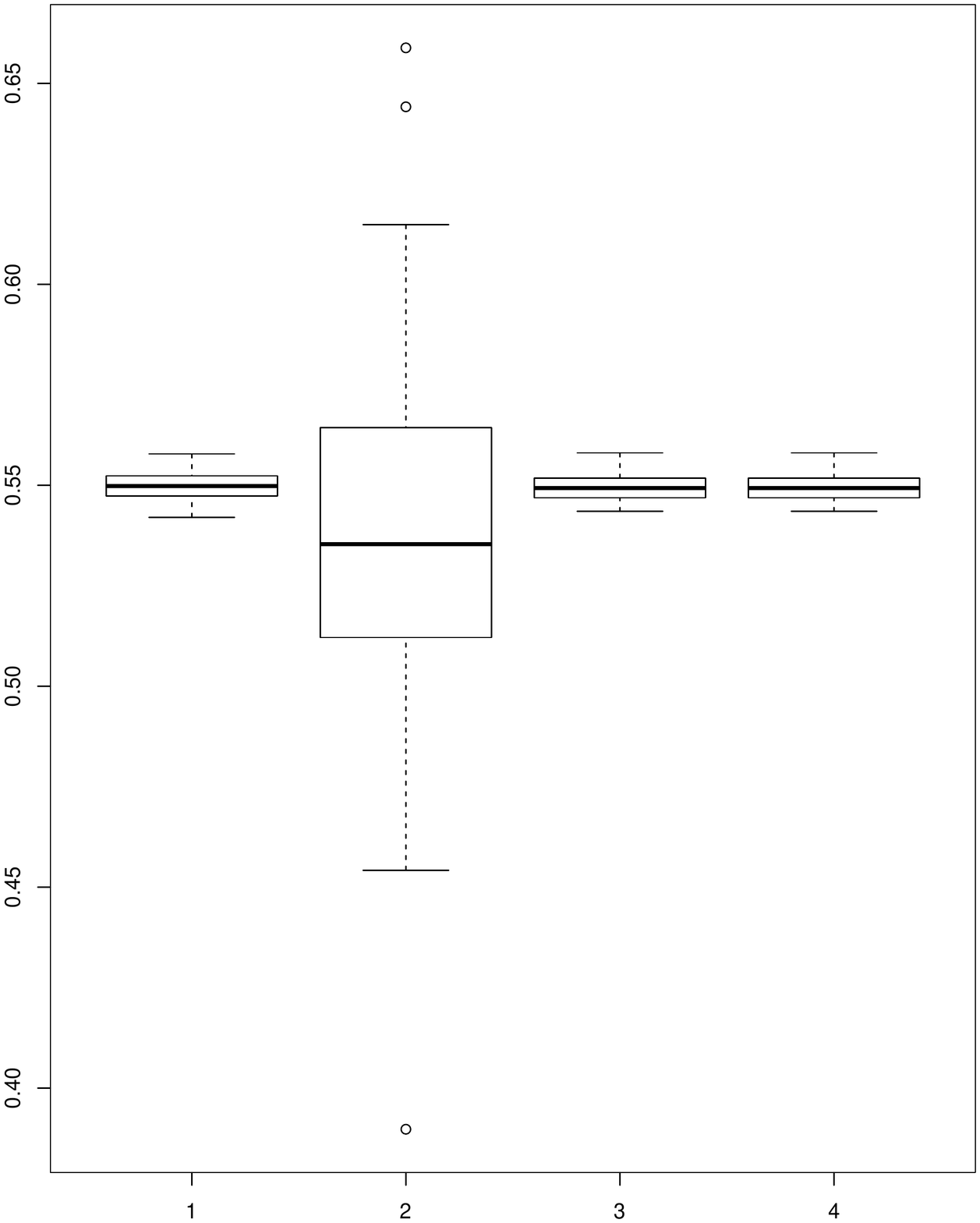}\quad\includegraphics[width=.22\linewidth]{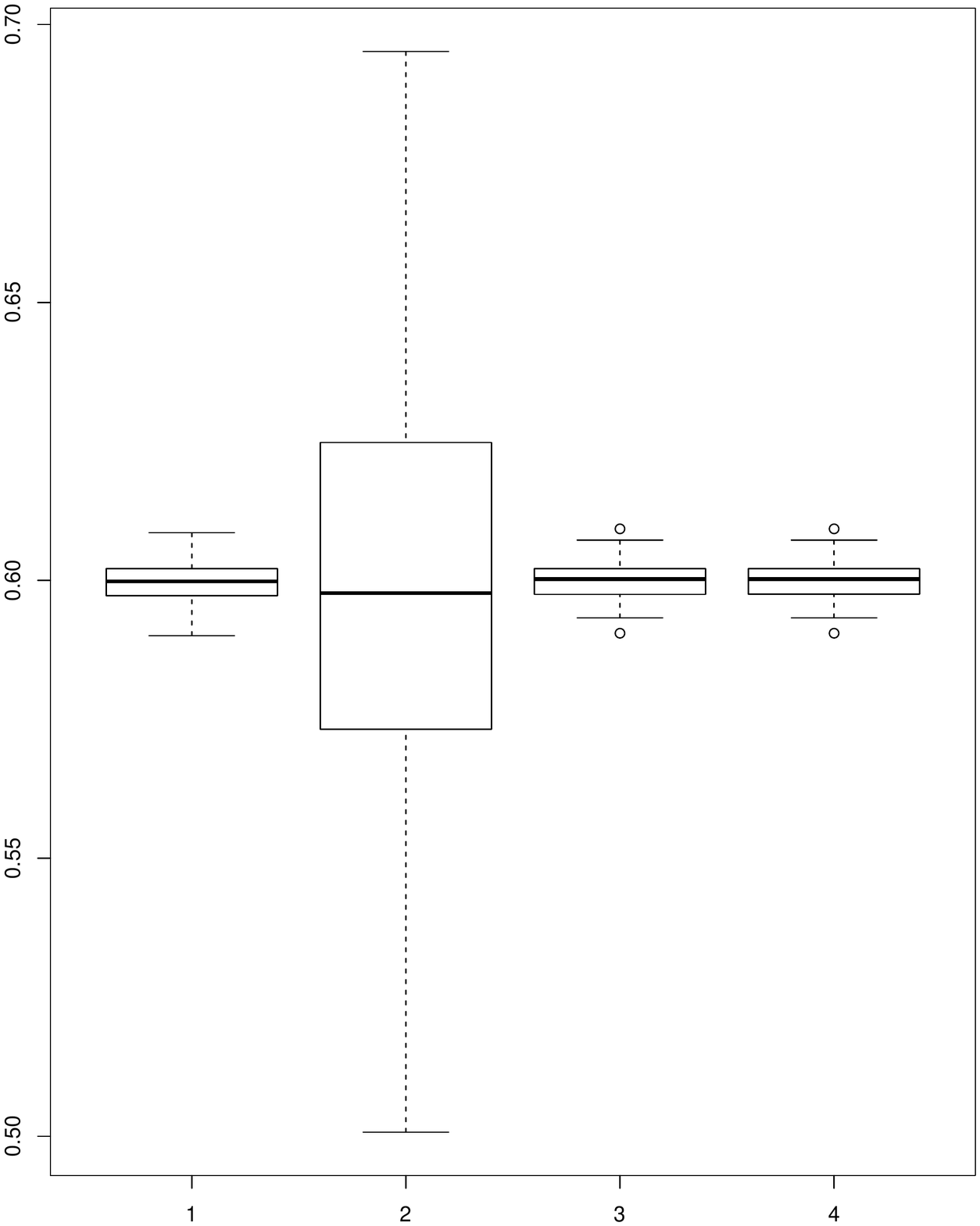}\quad \includegraphics[width=.22\linewidth]{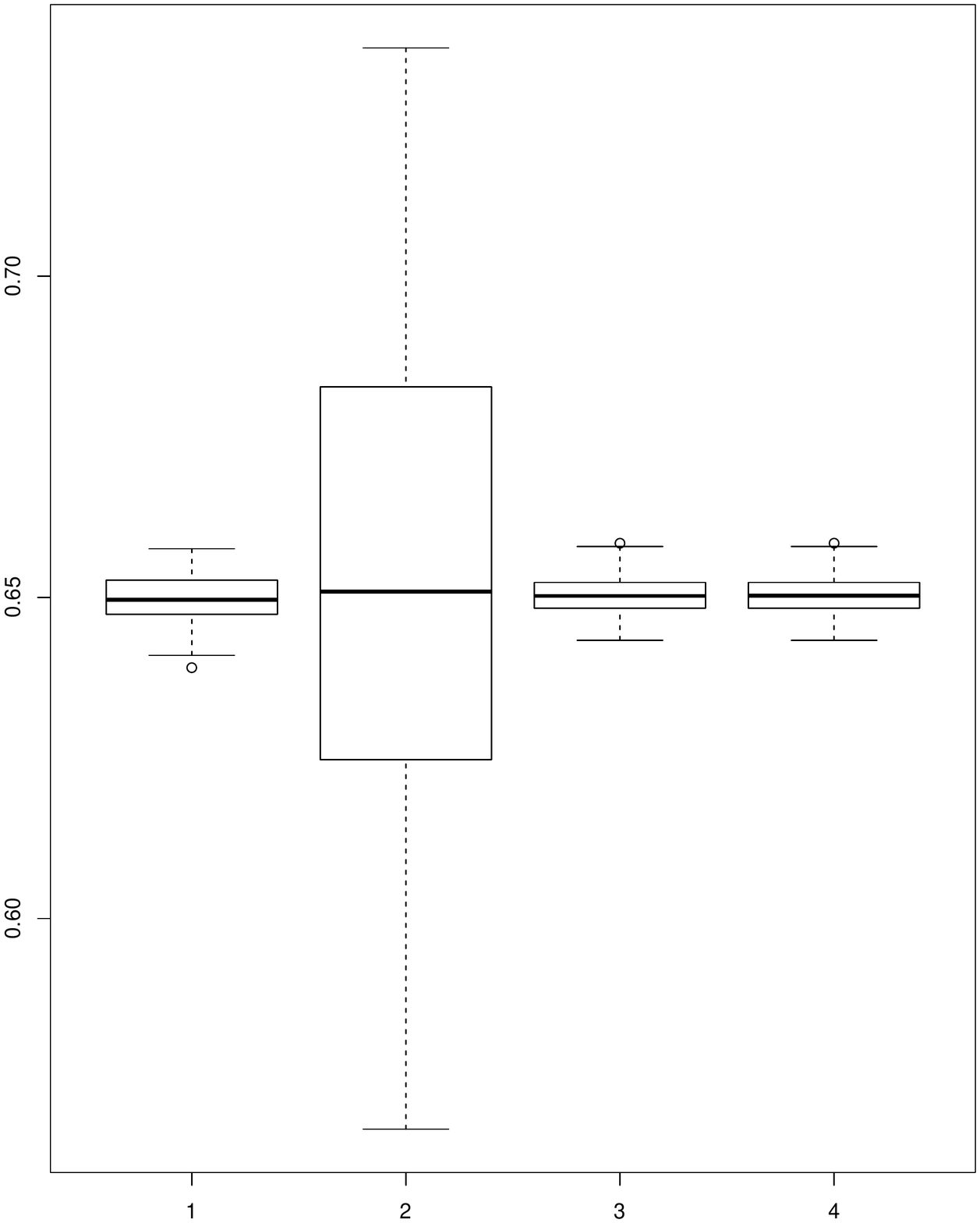}
\\[\baselineskip]
\includegraphics[width=.22\linewidth]{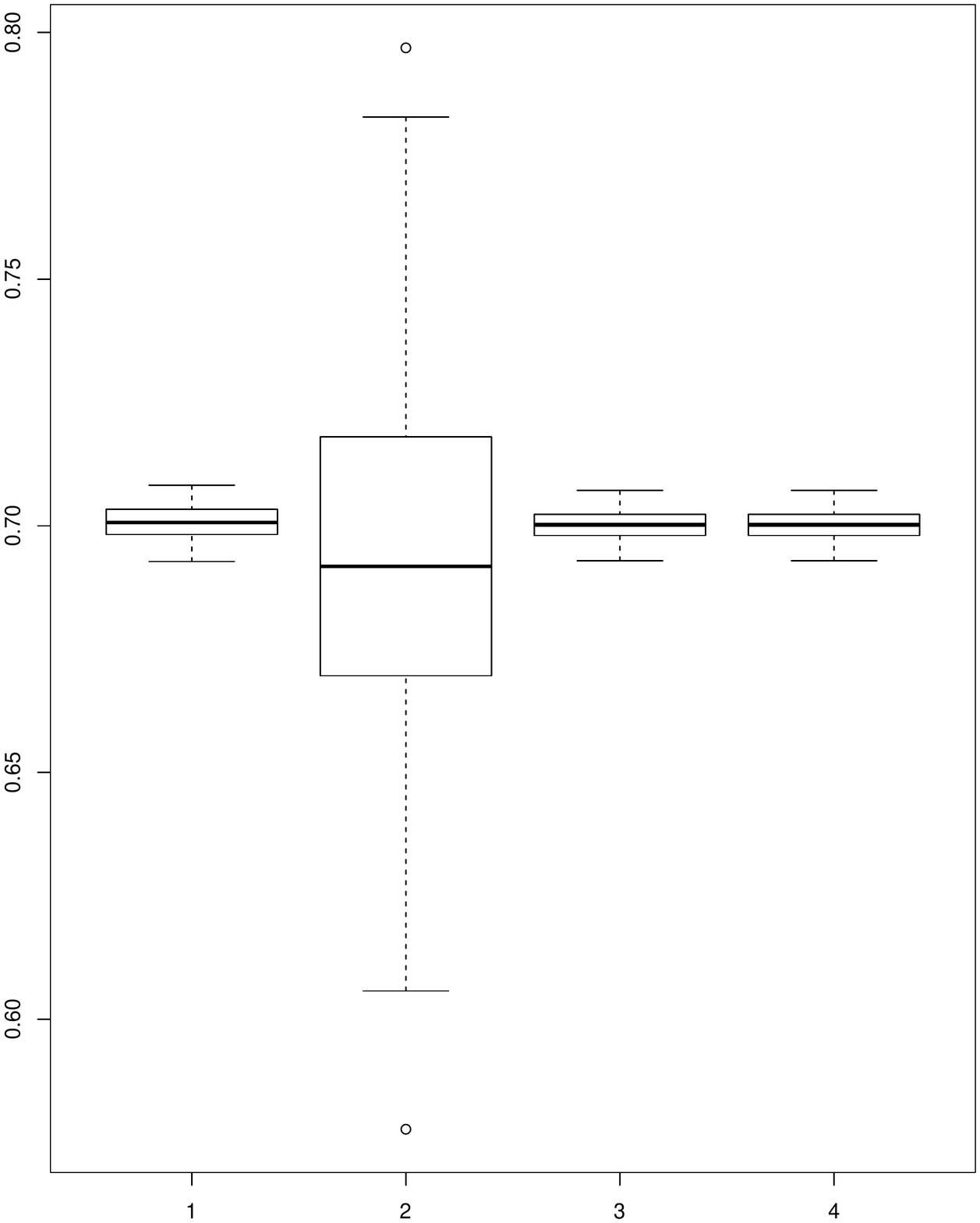}\quad\includegraphics[width=.22\linewidth]{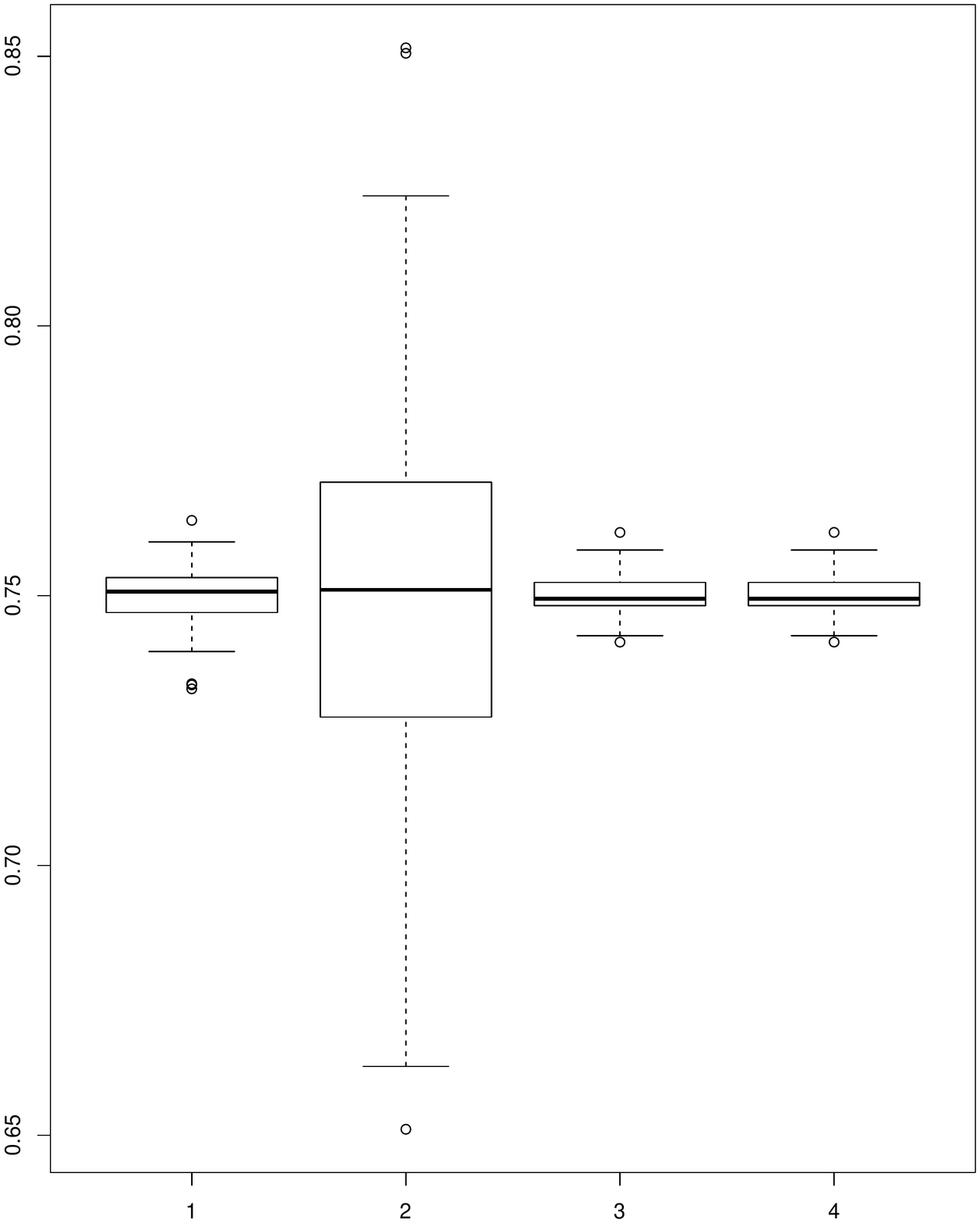}\quad \includegraphics[width=.22\linewidth]{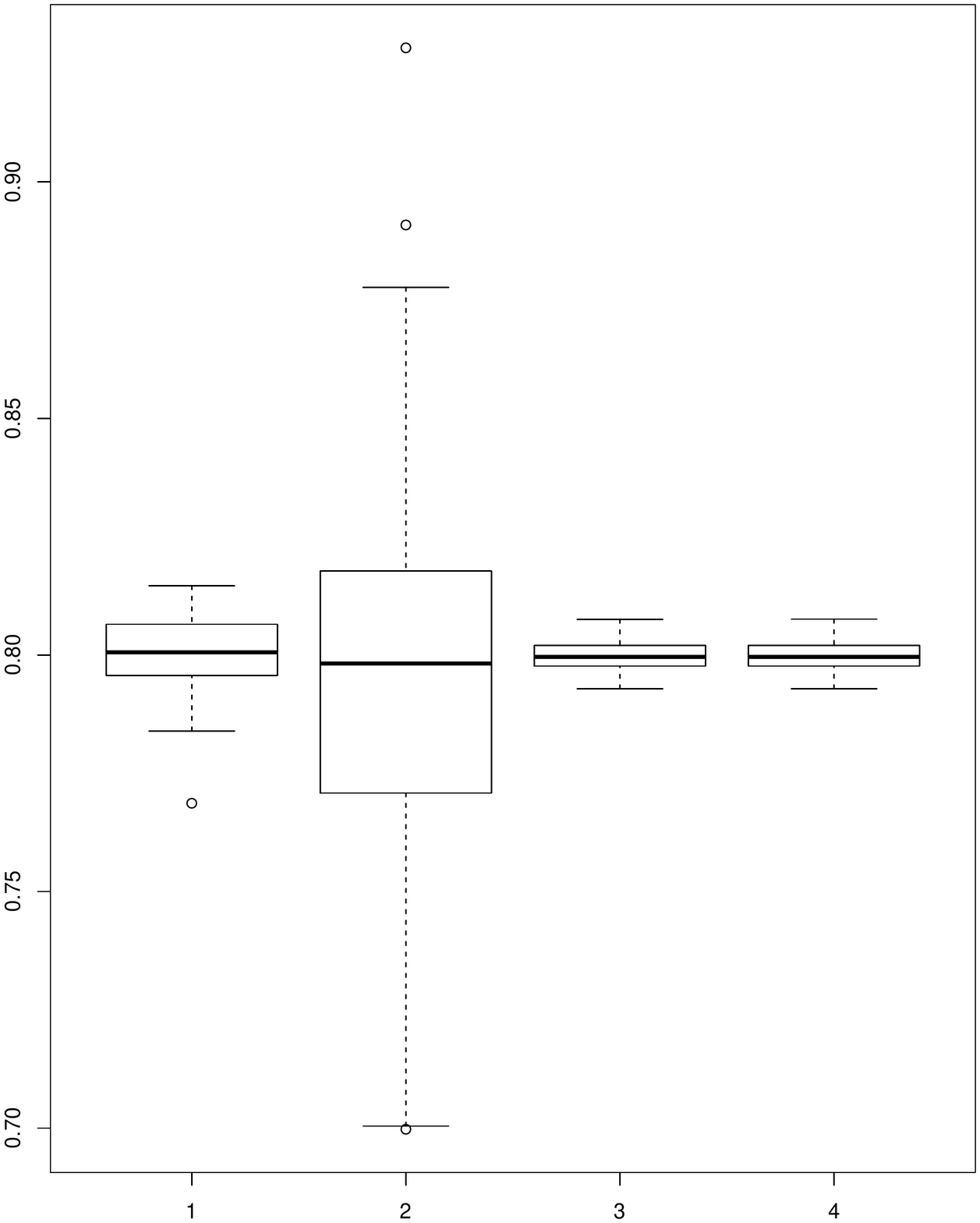}\quad\includegraphics[width=.22\linewidth]{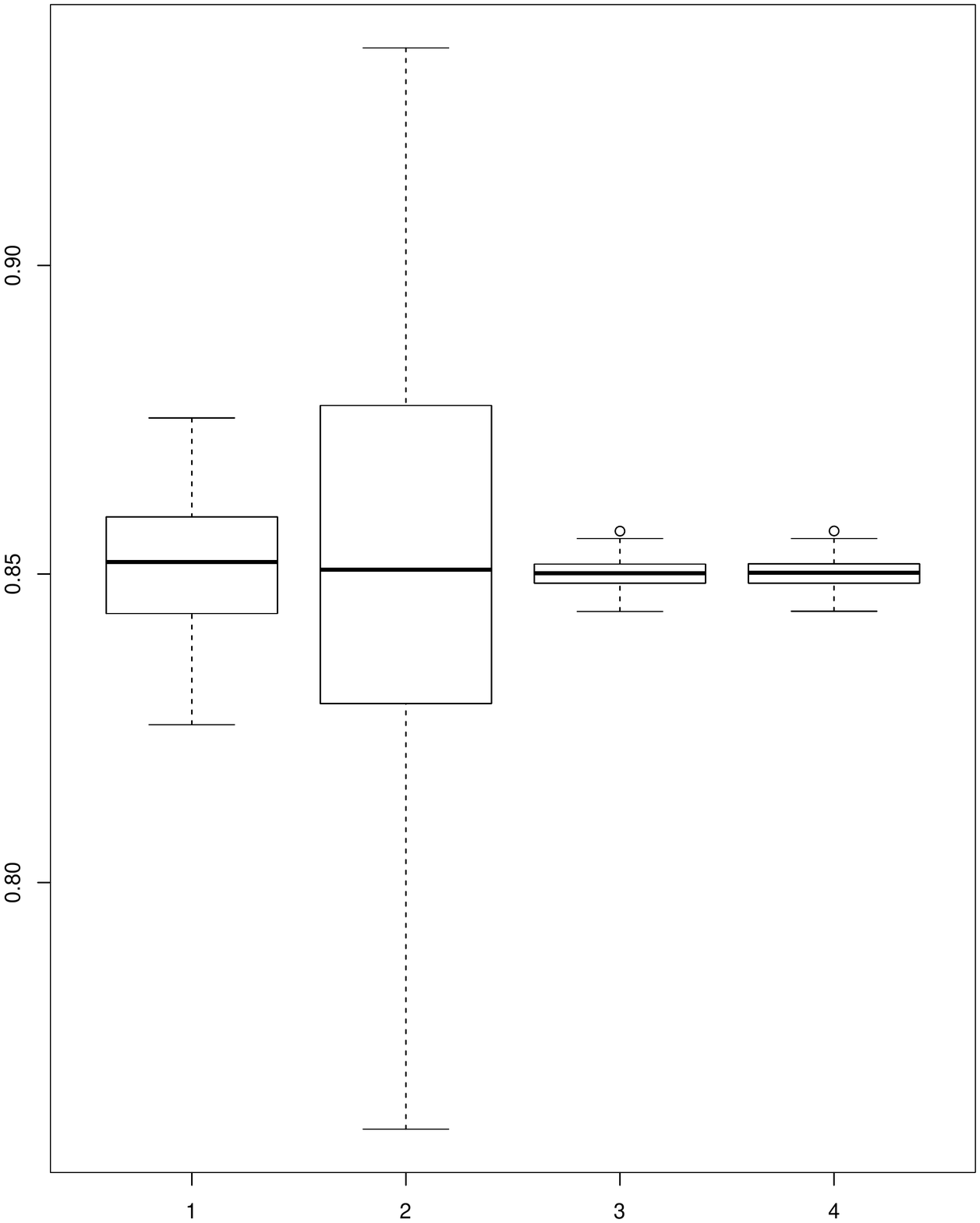}
\\[\baselineskip]
\includegraphics[width=.22\linewidth]{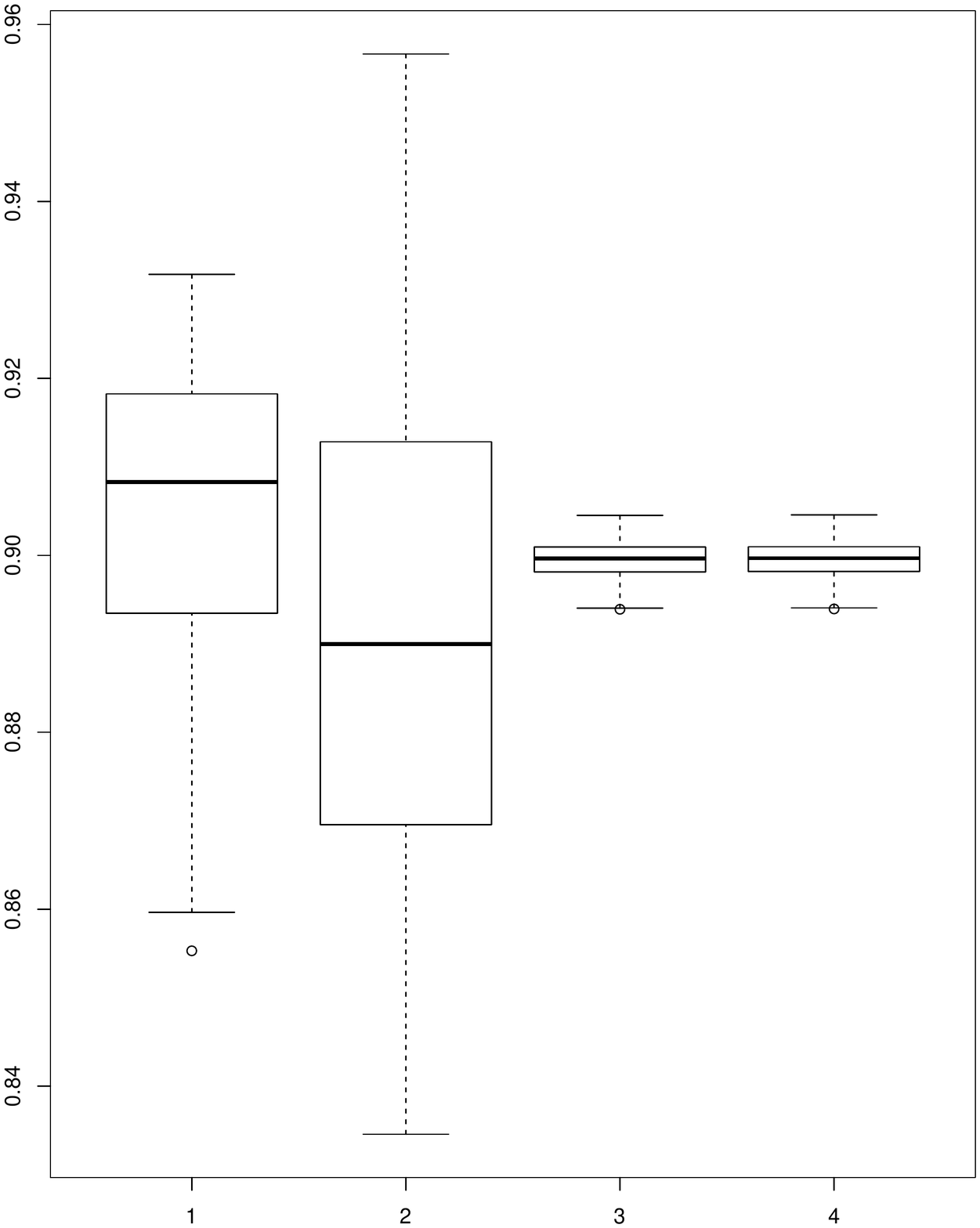}\quad \includegraphics[width=.22\linewidth]{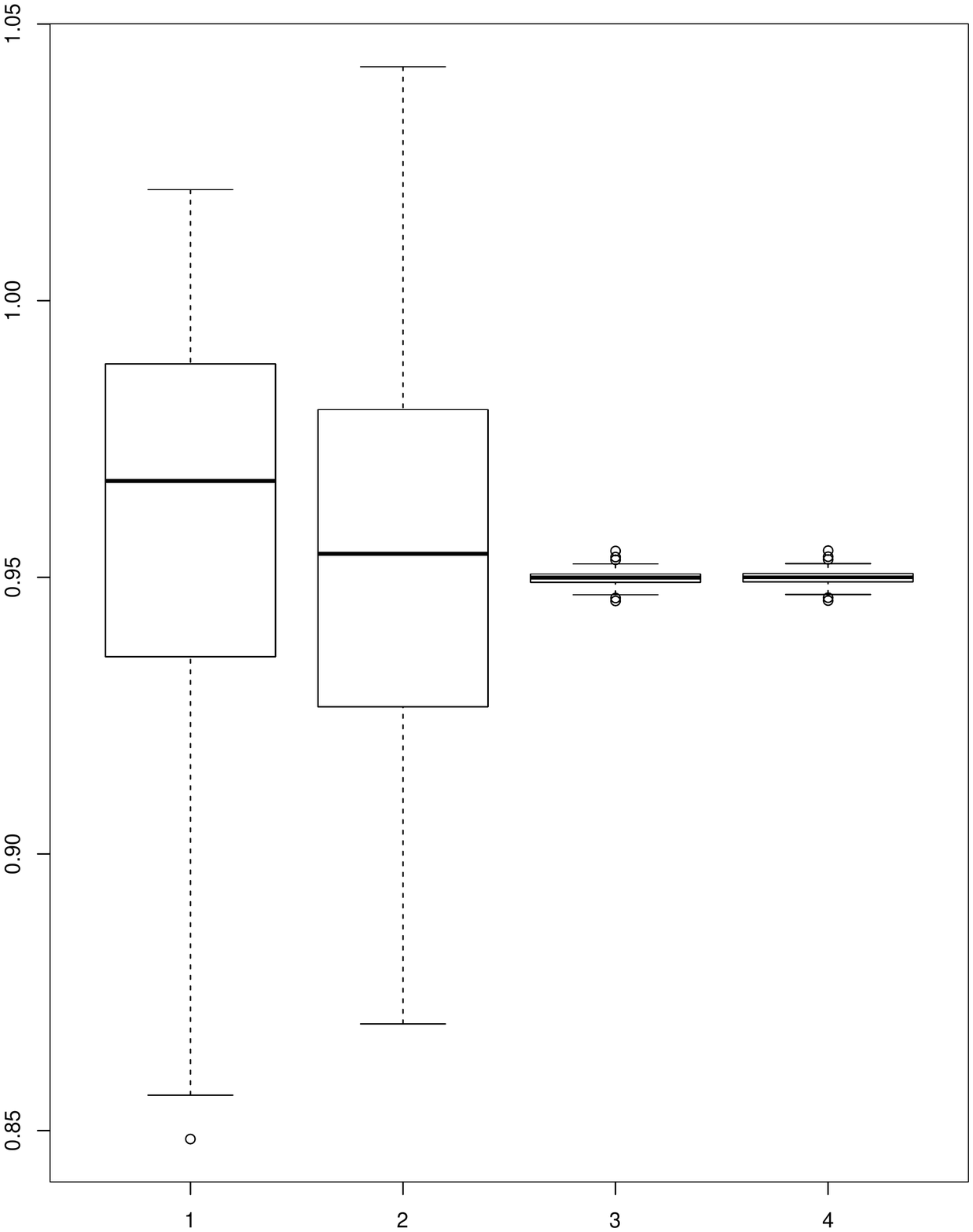}\quad\includegraphics[width=.22\linewidth]{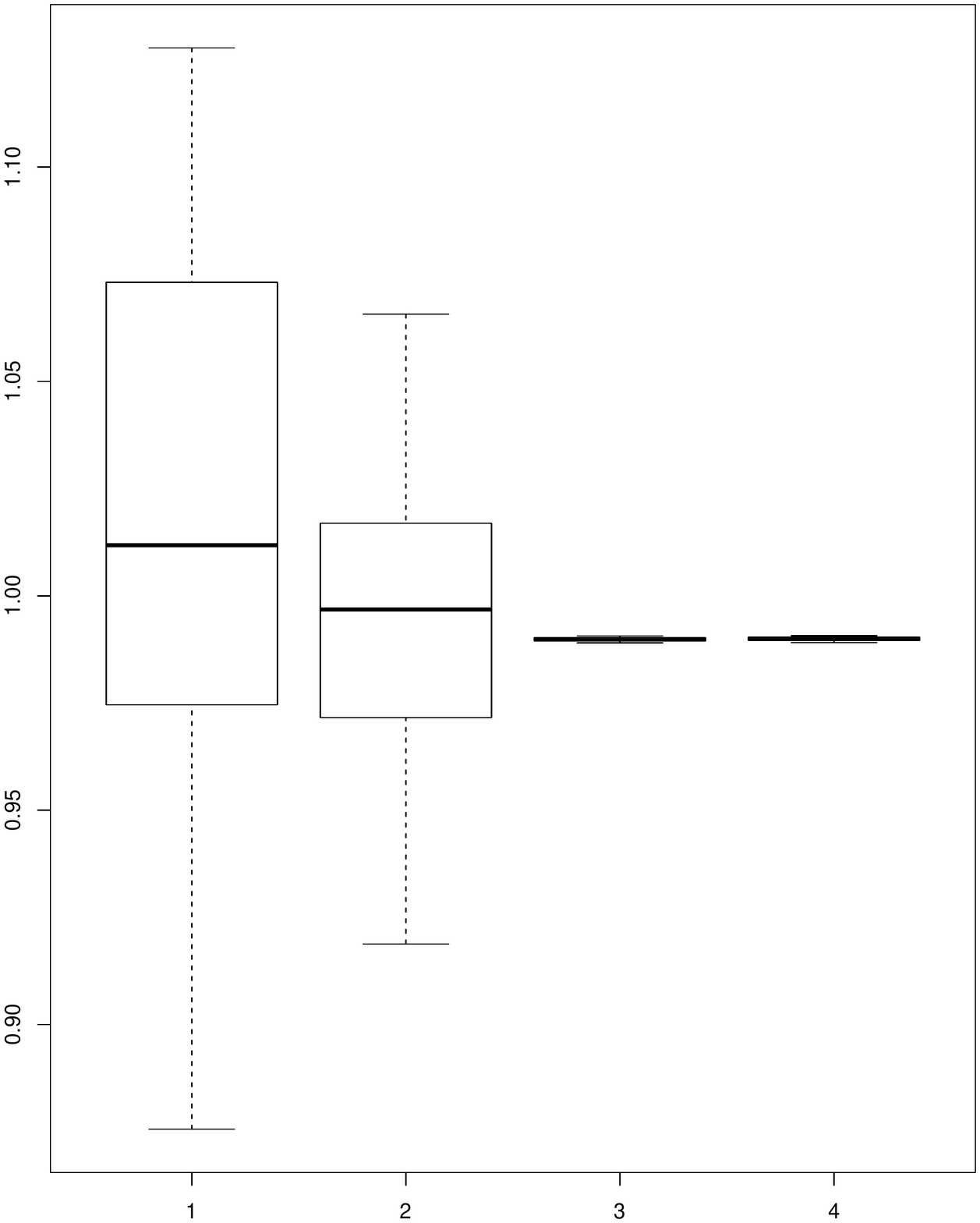}
\caption{Boxplots of  $\bar{H}_{2,\,N} (1,\,-1)$, $\bar{H}_{2,\,N}(1,\,-2,\,1)$, $\hat{H}_{2,\,N}(1,\,-2,\,1)$, $\tilde{H}_N$ for the values of $H$ listed above.}
\label{fig:bp}
\end{figure}

The Figures \ref{fig:paths} and \ref{fig:paths2} show the change in the limiting distribution, while the boxplots in Figure \ref{fig:bp} illustrate the changes in the speed of convergence indicated in the discussion for $\bar{H}_{2,\,N} (1,\,-1)$ and provide a comparison to the rates of convergence for the other three estimators.

\newpage

\section{Appendix}

The basic tools from the analysis on Wiener space are presented in this section. We will focus on some elementary facts about multiple stochastic integrals. We refer to \cite{N} for a complete review on the topic. 

Consider ${\mathcal{H}}$ a real separable infinite-dimensional Hilbert space
with its associated inner product ${\langle
.,.\rangle}_{\mathcal{H}}$, and $(B (\varphi),
\varphi\in{\mathcal{H}})$ an isonormal Gaussian process on a
probability space $(\Omega, {\mathfrak{F}}, \mathbb{P})$, which is a
centered Gaussian family of random variables such that
$\mathbf{E}\left( B(\varphi) B(\psi) \right) = {\langle\varphi,
\psi\rangle}_{{\mathcal{H}}}$ for every
$\varphi,\psi\in{\mathcal{H}}$. Denote by $I_{q}$ the $q$th multiple
stochastic integral with respect to $B$, which is an
isometry between the Hilbert space ${\mathcal{H}}^{\odot q}$
(symmetric tensor product) equipped with the scaled norm
$\frac{1}{\sqrt{q!}}\Vert\cdot\Vert_{{\mathcal{H}}^{\otimes q}}$ and
the Wiener chaos of order $q$, which is defined as the closed linear
span of the random variables $H_{q}(B(\varphi))$ where
$\varphi\in{\mathcal{H}},\;\Vert\varphi\Vert_{{\mathcal{H}}}=1$ and
$H_{q}$ is the Hermite polynomial of degree $q\geq 1$ defined
by:\begin{equation}\label{Hermite-poly}
H_{q}(x)=(-1)^{q} \exp \left( \frac{x^{2}}{2} \right) \frac{{\mathrm{d}}^{q}%
}{{\mathrm{d}x}^{q}}\left( \exp \left(
-\frac{x^{2}}{2}\right)\right),\;x\in \mathbb{R}.
\end{equation}The isometry of multiple integrals can be written as follows: for $p,\;q\geq
1$,\;$f\in{{\mathcal{H}}^{\otimes p}}$ and
$g\in{{\mathcal{H}}^{\otimes q}}$
\begin{equation} \mathbf{E}\Big(I_{p}(f) I_{q}(g) \Big)= \left\{
\begin{array}{rcl}\label{iso}
q! \langle \tilde{f},\tilde{g}
\rangle _{{\mathcal{H}}^{\otimes q}}&&\mbox{if}\;p=q,\\
\noalign{\vskip 2mm} 0 \quad\quad&&\mbox{otherwise}.
\end{array}\right.
\end{equation}It also holds that:
\begin{equation*}
I_{q}(f) = I_{q}\big( \tilde{f}\big),
\end{equation*}
where $\tilde{f} $ denotes the canonical symmetrization of $f$ and is defined by $$\tilde{f}%
(x_{1}, \ldots , x_{q}) =\frac{1}{q!}\sum_{\sigma\in\mathcal{S}_q}
f(x_{\sigma (1) },\ldots, x_{\sigma (q)}),$$ where the sum runs
over all permutations $\sigma$ of $\{1,\ldots,q\}$.

We have the following product formula: if  $f\in{{\mathcal{H}}^{\odot p}}$ and
$g\in{{\mathcal{H}}^{\odot q}}$, then 

\begin{eqnarray*}
I_{p}(f) I_{q}(g)&=& \sum_{r=0}^{p \wedge q} r! \binom{p}{r}\binom{q}{r}I_{p+q-2r}\left(f\tilde{\otimes}_{r}g\right)
\end{eqnarray*}
We need to recall the formula of contraction of elements of tensor products of Hilbert spaces. 
Consider  $\{e_{k}, k \geq 1 \} $ an orthonormal basis  of $\mathcal{H}$ and let $f\in{{\mathcal{H}}^{\otimes p}}$ and
$g\in{{\mathcal{H}}^{\otimes q}}$. 

For $r=1, \ldots, p\wedge q$,  the $r$th
contraction $f\otimes_{r}g$ is an element of
${\mathcal{H}}^{\otimes(p+q-2r)}$, which is defined by: 
\begin{eqnarray}\label{contra}
(f\otimes_{r} g)&=&\sum_{j_{1}, \ldots, j_{p}=1}^{\infty}<f, e_{k_{1}}\otimes e_{k_{2}}\otimes \ldots e_{k_{r}}>_{\mathcal{H}^{\otimes r}} \otimes <g, e_{k_{1}}\otimes e_{k_{2}}\otimes \ldots e_{k_{r}}>_{\mathcal{H}^{\otimes r}}. 
\end{eqnarray}
In the particular case when $\mathcal{H}=L^2(T)$, the $r$th
contraction $f\otimes_{r}g$ is the element of
${\mathcal{H}}^{\otimes(p+q-2r)}$ which is defined by {
\begin{eqnarray}\label{contra1}
& (f\otimes_{r} g) ( s_{1}, \ldots, s_{p-r}, t_{1}, \ldots, t_{q-r})  \notag \\
& =\int_{T ^{r} }\mathrm{d}u_{1}\ldots \mathrm{d}u_{r} f( s_{1},
\ldots, s_{p-r}, u_{1}, \ldots,u_{r})g(t_{1}, \ldots, t_{q-r},u_{1},
\ldots,u_{r})
\end{eqnarray}for every $f\in L^2(T^p)$, $g\in L^2(T^q)$
and $r=1,\ldots,p\wedge q$.
An important  property of  finite sums of multiple integrals is the hypercontractivity. Namely, if $F= \sum_{k=0} ^{n} I_{k}(f_{k}) $ with $f_{k}\in \mathcal{H} ^{\otimes k}$ then
\begin{equation}
\label{hyper}
\mathbf{E}\vert F \vert ^{p} \leq C_{p} \left( \mathbf{E}F ^{2} \right) ^{\frac{p}{2}}.
\end{equation}
for every $p\geq 2$.

We denote by $D$ the Malliavin derivative operator that acts on
cylindrical random variables of the form $F=g(B(\varphi
_{1}),\ldots,B(\varphi_{n}))$, where $n\geq 1$,
$g:\mathbb{R}^n\rightarrow\mathbb{R}$ is a smooth function with
compact support and $\varphi_{i} \in {{\mathcal{H}}}$. This
derivative is an element of $L^2(\Omega,{\mathcal{H}})$ and it is
defined as
\begin{equation*}
DF=\sum_{i=1}^{n}\frac{\partial g}{\partial x_{i}}(B(\varphi _{1}),
\ldots , B(\varphi_{n}))\varphi_{i}.
\end{equation*}

\end{document}